\documentclass[11pt]{article}
\input{amssym.def}
\input{amssym}
\newtheorem{theorem}{Theorem}

\newtheorem{definition}[theorem]{Definition}

\newtheorem{lemma}[theorem]{Lemma}

\newtheorem{proposition}[theorem]{Proposition}
\newtheorem{remark}[theorem]{Remark}

\def\uq{U_q(\gggg)}
\def\uqs{\ifmmode U_q(sl(n))\else $U_q(sl(n))$\fi}
\def\gggg{{\bf g}}

\def\C{{\Bbb C}}

\def\R{{\Bbb R}}
\def\K{{\Bbb K}}
\def\B{{\cal C}}
\def\lp{\Lambda_+(V)}
\def\lm{\Lambda_-(V)}
\def\lml{\Lambda_-^l(V)}
\def\lpm{\Lambda_{\pm}(V)}
\def\lpml{\Lambda_{\pm}^l(V)}

\def\la{{\lambda}}
\def\al{{\alpha}}
\def\vl{V_{\lambda}}
\def\cl{C_{\lambda}}
\def\shl
{s_{\lambda}}
\def\ot{\otimes}
\def\qm{q^{-1}}

\def\End{{\rm End\, }}

\def\Im{{\rm Im\, }}
\def\id{{\rm id\, }}
\def\tr{{\rm tr\, }}
\def\rk{{\rm rk\, }}

\def\Hom{{\rm Hom\,}}

\def\Mor{{\rm Mor}}
\def\Ob{{\rm Ob\,}}
\def\vv{V^{\ot 2}}
\def\SW{{\rm SW(V)}}
\def\Ru{\uq-{\rm Mod}}
\def\Rus{\uqs-{\rm Mod}}

\makeatletter \@addtoreset{equation}{section} \makeatother
\def\r#1{\mbox{(}\ref{#1}\mbox{)}}

\def\bea{\begin{eqnarray}}
\def\eea{\end{eqnarray}}
\def\beq{\begin{equation}}
\def\eeq{\end{equation}}
\def\be{\begin{equation}}
\def\ee{\end{equation}}

\def\mpi#1{\mbox{\boldmath  $\pi$}_{#1}}
\def\ev#1{\mathbf{ev}_{#1}}

\begin{document}

\input{amssym.def}
\input{amssym}

\title{\Large \bf Traces in braided categories}

\author{D.~Gurevich,
R.~Leclercq\\
{\footnotesize \it ISTV, Universit\'e de Valenciennes, 59304
Valenciennes, France}
\vspace{3mm}\\
P.~Saponov,\\
{\footnotesize\it Theory Department, Institute for High Energy
Physics, 142284 Protvino, Russia}}

\maketitle

\begin{abstract}
With any even  Hecke symmetry $R$ (that is a Hecke type solution
of the Yang-Baxter equation) we associate a quasitensor category.
We formulate a condition on $R$ implying that
 the constructed category is
rigid and  its commutativity isomorphisms $R_{U,V}$ are
natural in the sense of \cite{T}. We show that this
condition  leads to rescaling  the initial Hecke symmetry.
We suggest a new way of introducing traces as properly
normalized categorical morphisms $\End(V)\to \K$ and deduce
the corresponding normalization from categorical
dimensions.

\end{abstract}

\noindent {\bf AMS Mathematics Subject Classification, 1991 :}
18D10, 81R50


\section{Introduction}

The main purpose of this paper is to introduce some braided
categories arising from "nonquasiclassical" Hecke symmetries
constructed in \cite{G} and to suggest a way of defining
categorical traces \beq \tr=\tr_V: \End(V)\to\K ,\,\, V\in
\Ob(\B).   \label{tr} \eeq in a somewhat elementary way without
using any ribbon element (see below). Hereafter $\Ob(\B)$ stands
for the set of objects of the category $\B$ and $\K$ stands for
the basic field, always $\C$ or $\R$.

Let us precise from the very beginning that by braided categories
we mean monoidal tensor or quasitensor ones whose objects are
vector spaces with the usual tensor product and whose natural
associativity isomorphisms are identical (for the terminology,
cf. for example \cite{CP,T}). Thus, the structure of such a
category is principally determined by  commutativity isomorphisms
\beq
 R_{U,V}:U\ot V\to V\ot U,\qquad U, V\in \Ob(\B). \label{rm}
\eeq These isomorphisms are assumed to satisfy the so-called
Yang-Baxter (YB) equation. They play the role of the usual flip
which transposes the factors. We call the isomorphisms \r{rm} YB
operator, braiding, or quantum $R$-matrix.

Moreover, we assume them to be natural in the sense of \cite{T}.
This means that for any two categorical morphisms
 $f:\; U\rightarrow U'$ and $g:\;
V\rightarrow V'$ one has: \be (g\ot f)\circ
 R_{U,V} = R_{U',V'}\circ (f\ot g).\label{nat}
\ee We show that this condition leads to a normalization of
the commutativity isomorphisms different from that
(\ref{hec}) usually employed for Hecke type braidings. (Let
us note that even in the case of the quantum group
$U_q(sl(n))$ the YB operators coming from
 the universal R-matrix must be rescaled.)

In what follows we identify $\K$ and the $\End(I)$ where $I$ is
the identity object in the sense of \cite{CP}. Moreover,  we use
the notation $\End(V)=\Hom(V,V)$ for the objects of the category
which sometimes are called internal (endo)morphisms (cf.
\cite{DM}). By contrast, the morphisms of the category in question
will be called "categorical morphisms". For example, categorical
morphisms of the category $\gggg-{\rm Mod}$ of modules over a Lie
algebra $\gggg$ are elements of $\Hom(U,V)$, $U,
V\in\Ob(\gggg-{\rm Mod})$, commuting with the action of $\gggg$.

Initially, categories equipped with commutativity morphisms were
introduced by S.~MacLane \cite{ML}. However, he only considered
involutary YB operators ($R^2=\id$) called in what follows
"symmetries". A new wave of interest in braided categories (but
with noninvolutary braidings) arose in connection with integrable
system theory. More precisely, such categories appeared as those
of modules of Drinfeld-Jimbo quantum groups (QG) $\uq$ playing an
important role in this theory (we will denote this category
$\Ru$). Besides, the QG have found many other interesting
applications, in particular, in noncommutative geometry.

The problem we consider in the paper is a categorical
definition of the trace \r{tr}. What is a reasonable
generalization of the basic property of the classical trace
\beq \tr[X,\,Y]=0,\qquad \forall\, X, Y\in   \End(V)
\label{sled} \eeq where $[X,\,Y]$ is the commutator of two
endomorphisms\footnote{Sometimes this property is
attributed to the quantum trace (cf. \cite{H}). Equivalent
form of this property
 $\tr(g\, X\,g^{-1})=\tr X$ ($g\in
\End(V)$ is assumed to be invertible) appeared in \cite{T}.
However, all this is true if we assume $X,\, Y$ and $g$ to be
categorical morphisms. But in this case the conditions above
become
 meaningless. It is easy to see by taking a simple object $V$, since
it only admits scalar categorical morphisms and these conditions
degenerate.}?

If the braiding in a given braided category is involutary
and invertible by column in the sense of the formula
(\ref{R-Q}) then there exists a natural generalization of
the above commutator such that relation \r{sled} is still
valid for this commutator and an appropriate trace. For
example it is so for a super-commutator and super-trace
defined in a super-category. Other examples can be found in
\cite{G}\footnote{In this case the trace can be treated in
termes of "R-cyclic cohomology" which can be naturally
defined via the operator $R$ (in the spirit of super-cyclic
 cohomology). However, apparently there does not exist any reasonable
similar  treatment  of the trace in the categories under
consideration.}.

Let us note that the definition of the trace \r{tr} in a category
equipped with a symmetry makes use of this symmetry and of the
identification $\End(V)=V\ot V^*$ where $ V^*$ is an object dual
to $V$ (all categories in question are assumed to be rigid, i.e.
for any object of the category its dual is also an object). Up to
our knowledge for the first time such a trace has been introduced
in \cite{DP}.

Once such a trace is defined one can introduce {\it a categorical
dimension} $\dim(V)$ of an object $V$ by setting $\dim(V) =
\tr(\rm id)$. It is easy to see that \beq {\dim} (U\oplus
V)={\dim}(U)+{\dim}(V),\,\,{\dim}(U\otimes V)={\dim}(U){\dim}(V),
\label{a-m} \eeq i.e. the dimension can be considered as an
additive and multiplicative (a-m) functional on the objects of
the category.

However, a direct application of the approach of \cite{DP} to the
category $\Ru$ gives "dimension" which is not an a-m functional.
This is the reason why one has to introduce a correction in the
definition of "dimensions" turning them into an a-m functional.
The correction is connected to the so-called ribbon element in the
corresponding QG. The image of this element is sometimes called
"twist" (cf. \cite{T}). The categories possessing a braiding and a
twist are called ribbon. The trace in such categories is defined
via some combination of the ribbon element and another element $u$
due to Drinfeld (cf. \cite{CP}) and this leads to the dimension
which is an a-m functional.

Some generalization of this construction of the trace has been
suggested in \cite{BW}. Suggested in that paper is a way to define
the notion of dimension in some categories without any braiding.
Such a category is introduced as that of $H$-modules where $H$ is
a Hopf algebra. Instead of a universal quantum $R$-matrix (i.e.
the element  giving rise to a braiding) the authors of \cite{BW}
use another element allowing to identify any object $V$ with its
second dual $V^{**}$.

In the present paper we restrict ourselves to categories equipped
with braidings and suggest another somewhat elementary way of
introducing traces.

Before discussing this way let us introduce some notations and
definitions.

A Yang-Baxter operator \beq R:\vv\to\vv  \label{rmat} \eeq will
be called {\it a Hecke symmetry} if in addition to the
Yang-Baxter equation \beq
R_{12}\,R_{23}\,R_{12}=R_{23}\,R_{12}\,R_{23} \label{YB} \eeq it
also satisfies the relation \beq (q\,\id-R)(q^{-1}\,\id+R)=0,\quad
q\in \K. \label{hec} \eeq Here the standard tensor notations have
been used for the equations and operators in tensor product of
spaces. Equation (\ref{YB}) is written in $V^{\ot 3}$ and one
assumes
$$
R_{12} = R\otimes {\rm id}, \quad R_{23} = {\rm id}\otimes R,
$$
where ${\rm id}$ is the identity operator on $V$.

The nonzero parameter $q$ is taken to be generic, which means  it
is not a root of unity: $\forall k\in{\Bbb N},\;k\not=1,\;q^k\not
=1$. As a consequence none of the so-called {\it $q$-numbers}
$k_q$ is equal to zero \be k_q\equiv \frac{q^k -
q^{-k}}{q-q^{-1}} \not=0\quad\forall k\in\Bbb N.\label{q-num} \ee

Let us assign to such a Hecke symmetry "symmetric" $\lp$ and
"skew-symmetric" $\lm$ algebras by \be
\lp=T(V)/\{\Im(q\,\id-R)\},\quad \lm=T(V)/\{\Im(\qm\id+R)\}.
\label{Sym-alg} \ee Here $T(V)$ stands for the free tensor algebra
of the space $V$ and $\{I\}$ denotes the ideal generated  by a
subset $I$ in a given algebra. Let $\lpml$ be a homogeneous
component of the algebra $\lpm$ of the degree $l$. We call a Hecke
symmetry {\it even} if there exists an integer $p$ such that the
component $\lml$ is trivial at $l>p$  and it is one-dimensional at
$l=p$. The integer $p$ will be called a rank of $V$ and be denoted
$\rk(V)$.

In the latter 80's one of the authors (D.G.) constructed
examples of Hecke type braidings which differ drastically
from those related to the QG (see \cite{G} and references
therein). Namely, it was shown that there exist a lot of
Hecke symmetries $R$ such that $n=\dim(V)>p=\rk(V)$. Here
by $\dim(V)$ we mean the usual dimension of the space $V$.
Note that these Hecke symmetries are not deformations of
the usual flip and we call them
nonquasiclassical\footnote{Some of them have been
independently introduced in \cite{DL}.}. To the contrary,
the Hecke symmetry coming from the QG $\uqs$ is a
deformation of the flip and we call it quasiclassical.

In this paper we introduce a braided category generated by a
vector space $V$ equipped with a Hecke symmetry such that its
braidings are natural. We call it the Schur-Weyl (SW) category
and denote as $\SW$. Let us remark that similar categories were
considered in some papers (cf. \cite{KW}, \cite{H}, \cite{B}).
However, if the authors of \cite{KW} "reconstruct" an existing
category we have no category at the very beginning and should
first construct it.

We construct the category directly by giving the list of objects
and categorical morphisms without using any RTT algebra
habitually employed for that (\cite{H}, \cite{B}). This leads to
different defining morphisms of the category and finally, to the
condition (\ref{fix-MN}) ensuring naturality of the braidings
$R_{U,V}$. (In a separate publication one of us (R.L.) shows that
in the class of nonquasiclassical Hecke symmetries from \cite{G}
there exists a big subclass of those satisfying this condition.)

Now let us describe the category under consideration. Any object
of such a category is a direct sum of simple ones and in this
sense the category $\SW$ is spanned by simple (basic) objects.
These basic objects are labelled by  partitions (or what is the
same by Young diagrams) \beq \lambda =(\lambda _{1},\lambda
_{2},...,\lambda _{p-1}),\quad\lambda _{i+1}\le \lambda _{i},
\label{part} \eeq $\lambda_i$ being nonnegative integers. In the
sequel we will use the notation $\lambda\vdash k$ for each
$\lambda$ being a partition of the integer $k$ that is
$\sum_i\lambda_i = k$. The number of nonzero components of a
partition $\lambda$ is called its {\it height} and will be denoted
$\ell (\lambda)$. For an object corresponding to a diagram $\la$
we introduce the notation $\vl$.

Let us stress that the tensor product of two basic objects is a
direct sum of basic ones. This naturally leads us to the notion of
Grothendieck semiring (and hence ring) of the category which
turns out to be the same as for the categories of $sl(n)$- (or
$\uqs$)-modules. The only difference is that the role of the
classical dimension $n=\dim(V)$ is played by the rank $p=\rk(V)$.
Thus, we can consider the dimensions of objects of the category
as a functional on the Grothendieck semiring. Since it is an a-m
functional we come to the problem of describing the whole family
of a-m functionals and select those related to  traces.

By using the standard technique of symmetric functions (cf.
\cite{M}) we show in Proposition~\ref{AM-func}  that such a
functional $f$ is determined by its values on objects
$V_{(1^k)}$, $k=1,2,...,p-1$, and its value on a basic object
$V_\lambda$ is equal to
$$
f(V_\lambda) = \shl(\al_1,\al_2,..., \al_p)
$$
where $\shl$ is the Schur function (polynomial) in $p$ variables
corresponding to the partition $\la$. The numbers $(-\al_i)$ are
the roots of the polynomial
$$
\phi(t) = t^p + f_1t^{p-1} +\dots + f_{p-1} t +1,\qquad
f_k=f(V_{(1^k)}).
$$

(Note, that the classical (usual) dimension is also an a-m
functional. Therefore, the above result allows one to calculate
the classical dimension of spaces $V_\lambda$ provided that the
Poincar\'e series $P_-(t)$ (see Remark~\ref{concl}) is known (for
involutary symmetries such a calculation has been done in
\cite{GM}). Remark, that if a Hecke symmetry $R(q)$ (\ref{rmat})
is a deformation of an involutary one $R(1)$ then the usual
dimension of corresponding spaces $\vl$ will be the same in both
cases. Indeed, being an integer number such a dimension is stable
under deformation, while the categorical dimension is a function
in $q$ which becomes integer at $q=1$.)

Having classified all a-m functionals one can put the question:
which a-m functionals come from SW categories? More precisely, for
which a-m functional on a given Grothendieck semiring (with a
fixed $p$) there exists a SW category such that the categorical
dimension of its objects coincides with this functional?

Now we are going back to the problem of defining traces on objects
$\End(U)=U\ot U^*$ (for the sake of concreteness we consider algebras
of left endomorphisms). Let us  assume  that trace is defined on
$\End(U\oplus V)=(U\oplus V)\ot (U^*\oplus  V^*)$
via
$$\tr\vert_{\End(U)\oplus \End(V)}=\tr\vert_{\End(U)}+\tr\vert_{\End(V)}$$
being extended to $U\ot V^*\oplus V\ot U^*$ by 0.
Then traces are completely defined by their values on
simple objects. If $V$ is such an object then the trivial
component in the product $V\ot V^*$ is unique. So, the
morphism \r{tr} for such $V$ (being nothing but a multiple
of the projection on the trivial component) is uniquely
defined up to a nontrivial factor.

As a trace we  take this properly normalized morphism.  Our
"proper" normalization is fixed by the requirement that the
corresponding dimension should be an a-m functional on the
corresponding Grothendieck semiring. Thus, this normalization is
a collective phenomenon.

Remark that our method to define  traces as properly
normalized categorical morphisms is valid, in principle,
for any  category for which the Grothendieck semiring is
well-defined. This method does not make  use of either
ribbon element or twist. It is sometimes useful to employ
twists in order to calculate categorical dimensions (as it
is done in \cite{H}) but it is not reasonable to introduce
it in the definition of  traces. Also, it turns out that
this method is useful for studying K-theory of some
algebras related to SW categories and computing
noncommutative index for them. This application will be
presented elsewhere.

As for  the category under consideration a computation shows that
the corresponding dimensions of basic objects are
 \beq {\dim}_q(\vl)=s_{\lambda
}(q^{p-1},q^{p-3},...,q^{-(p-1)}). \label{dim} \eeq In a little
bit different form the formula (\ref{dim}) was also given in
\cite{KW} and \cite{H}. (Also remark that  in the case of the
category $\Rus$ the above formula is  equivalent to  that from
\cite{Go}.)

However, behind this nondetailed information about the
objects encoded in their dimensions different and rich
structures are hidden. For example, traces (\ref{tr}) in
the categories in question are defined via some matrix $C$
entering the formula for the trace on the basic space $V$
and extended in a proper way onto the whole category (the
extension of $C$ matrix to the space $\vl$ will be denoted
$\cl$). It is worth mentioning that the matrices $\cl$ are
essentially different in quasiclassical and
nonquasiclassical cases while the categorical dimensions of
objects which are nothing but the properly normalized usual
traces of $\cl$ depend only on $\la$ and $p=\rk(V)$ via the
formula \r{dim}.

The paper is organized as follows. In Section~2 we
construct all the necessary elements of a Schur-Weyl
category: the class of objects, categorical morphisms, the
quasitensor and rigid structure. At the end of Section~2 we
discuss the role of the condition (\ref{fix-MN}) for
applications. In Section~3 we describe the full set of a-m
functionals on the objects of SW(V), define the trace in
End(U), $\forall\, U\in {\rm Ob(SW(V))}$, and calculate
$\dim_q U$ as the value of trace on identity morphism in
End(U).
\bigskip

{\bf Acknowledgement}\ Two of the authors (D.G. and P.S) would
like to thank Max-Planck-Institut f\"ur Mathematik (Bonn) where
the paper was completed for warm hospitality.

\section{Schur-Weyl category generated by a Hecke symmetry}

In this section we describe the construction of Schur-Weyl
category SW(V) (see Introduction) generated by a finite
dimensional vector space $V$ equipped with a Hecke symmetry. We
supply the category with attributes of that \Rus{}. In particular,
all categories SW(V) are $\Bbb K$-linear rigid quasitensor ones
(for precise definitions see \cite{CP,MacL}). As was shortly
outlined in Introduction, the basic objects of SW(V) are some
linear subspaces $V_\lambda$ of the tensor algebra $T(V)$, while
the categorical morphisms will be specific linear mappings of the
objects (see below). Consequently, a $\Bbb K$-linear structure of
such a category is obvious: as the direct sum of objects we will
take the usual direct sum of linear spaces and the trivial (zero)
vector space will be the null object of the category. So, we
should only define the structure of monoidal category and fix
braiding isomorphisms converting the monoidal structure into the
quasitensor (or braided) one. Besides, to have a rigid category,
we must ensure that the dual space $V^*$ of an arbitrary $V\in
{\rm Ob(SW(V))}$ is also an object of the category. Now we pass
to the explicit construction of the mentioned components of the
category SW(V) taking into account condition (\ref{nat}).

\subsection{\it Class of objects}

Let us fix some finite dimensional vector space $V$, $\dim V = n$,
and choose a basis $\{e_i\}$ in this space. Consider a linear
operator $R$ (\ref{rmat}) whose action on the basis elements
reads:
$$
R(e_i\otimes e_j) = e_r\otimes e_s\,R^{rs}_{\,\, ij},
$$
where the summation over the repeated indices is understood. In
what follows we assume $R$ to be a Hecke symmetry with a generic
$q\in \Bbb K$.

Besides, we will suppose $R$ to be an even symmetry of a rank
$p\le n = {\dim V}$. This means that the Poincar\'e series
$P_-(t)$ related to the "skew-symmetric" algebra $\wedge_-(V)$
(\ref{Sym-alg}) is a monic polynomial of $p$-th degree, that is a
polynomial in $t$ whose leading coefficient is equal to~1.

Throughout the paper we will use the compact notation for the
relations including matrices and vector spaces. In such notations
the indices will stand for a number of space rather than for a
particular  matrix or vector component. A basis vector $e_i$ in
the $k$-th matrix space is denoted as $e_{\langle k|}$, while for
the basis vector of dual space we will write $e^{|k\rangle}$. The
same principle is applied to components of arbitrary tensors.

For example, the above relation will look like\footnote{The symbol
$:=$ means "by definition".}
$$
R(e_{\langle 1|}\otimes e_{\langle 2|}) := e_{\langle 1|}\otimes
e_{\langle 2|}\,R^{|12\rangle}_{\;\;\langle 12|}:= e_{\langle
1|}\otimes e_{\langle 2|}\,R_{12}.
$$
In this notation the summation over the repeated indices is
represented  by an expression which contains the same upper and
lower case indices {\it with properly oriented brackets} as shown
in the examples below:
$$
e_{\langle 1|}\ot e^{|1\rangle}:= \sum_i e_i\ot e^i,\quad
e_{\langle1|}\cdot T_1 := e_{\langle1|}\,T^{|1\rangle}_{\, \langle
1|}:= \sum_j e_jT^j_{\,i},\quad M_1R_{12}N_2 :=
\sum_{a,b}M^{i_1}_{\;a}R^{ai_2}_{\,\,j_1b}\,N^b_{\,j_2}
$$
but
$$
e^{|1\rangle}\ot e_{\langle 1|} := e^{i}\ot e_{j},\quad e_{\langle
1|}\ot e^{|2\rangle} := e_i\ot e^j,\quad T_1\cdot e_{\langle1|} :=
T^i_{\,j}\,e_{k}
$$
and so on.

Now we use the fact, that in each homogeneous component $V^{\ot
m}$ of the tensor algebra $T(V)$ one can realize the so called
{\it local representation} of the Hecke algebra $H_m$
\cite{DipJ}  via a given Hecke symmetry\footnote{For the review
on Hecke algebra the reader is referred to the recent work
\cite{PyOg}}.

The Hecke algebra $H_m$ is generated by the unit element
$\id_{\!\! H}$ and $m-1$ generators $\sigma_k$ subject to the
following relations:
$$
\left.
\begin{array}{l}
\sigma_i\sigma_{i+1}\sigma_i = \sigma_{i+1}\sigma_i\sigma_{i+1}\\
\sigma_i\sigma_j = \sigma_j\sigma_i\qquad {\rm if\ }|i-j|\ge 2\\
(\sigma_i-q\,\id_{\!\!H})(\sigma_i+q^{-1}\,\id_{\!\!H}) = 0\\
\end{array}
\right\}\;\; i=1,2,\dots m-1.
$$
The local representation of $H_m$ in $V^{\ot m}$ is of the form:
\be \sigma_i\rightarrow \rho_R(\sigma_i)\equiv R_{ii+1} = {\rm
id}_{i-1}\otimes R\otimes {\rm id}_{m-i-1}\in {\rm End}(V^{\ot
m}).\label{Rep} \ee

In the further construction of SW(V) category the central role
belongs to the fact that the Hecke algebra $H_m$ for a generic $q$
is semisimple and can be decomposed into a direct sum of ideals.
Moreover, since the Hecke algebra $H_m$ at a generic $q$ is
isomorphic to the group algebra ${\Bbb K}[S_m]$ of the  $m$-th
order permutation group $S_m$, its primitive idempotents
generating the ideals in mentioned decomposition can be put into
one-to-one correspondence with the set of all standard Young
tableaux connected with each possible partition $\lambda $ of the
integer $m$. Speaking more explicitly, given the algebra $H_m$ and
a partition $\lambda\vdash m$, one can construct some polynomials
$Y_{ii}^{\lambda}$ in generators $\{\sigma_i\}$ which turn out to
be the primitive idempotents of Hecke algebra $H_m$ (for detailed
description of such a construction see review \cite{PyOg}): \be
\id_{\!\!H} = \sum_{\lambda\vdash
m}\sum_{i=1}^{d_{\lambda}}Y_{ii}^{\lambda},\label{id-dec} \ee

\be Y_{ii}^{\lambda}\,Y_{jj}^{\mu} =
\delta_{ij}\delta^{\lambda\mu}\,Y_{ii}^{\lambda}.\label{Y-ortog}
\ee In the above formulas the index $i$ runs from 1 to
$d_{\lambda}\equiv \dim\lambda$, that is to the number of all
standard Young tableaux, corresponding to the Young diagram of the
partition $\lambda\vdash m$. Recall, that a Young tableau is
called {\it standard} if it is filled with successive integers
from $1$  to $m$ in such a way that the integers increase  from
left to right in each row and from top to bottom in each column.
The standard Young tableaux (and hence the primitive idempotents
$Y^\lambda_{ii}$) may be lexicographically ordered in many ways
and for definiteness we fix an order in which the first idempotent
$Y^\lambda_{11}$ corresponds to the Young tableau filled by
integers consequently increasing by $1$ when going down in each
column. Here is an example for the Young tableau corresponding to
the partition $\lambda = (3,2^2,1)$: \be Y_{11}^{\lambda}
\longleftrightarrow\;
\begin{array}{|c|c|c|}\hline
1&5&8\\ \hline 2&6&\multicolumn{1}{c}{}\\ \cline{1-2}
3&7&\multicolumn{1}{c}{}\\ \cline{1-2} 4&\multicolumn{2}{c}{}\\
\cline{1-1}
\end{array}
\label{Example} \ee

A primitive idempotent corresponding to the partition
$\lambda=(1^m)$ will be called the antisymmetrizer and be denoted
as $A^{(m)}$. There is one important circumstance about the
images of $A^{(m)}$ in End$(V^{\ot m})$ with respect to
representation (\ref{Rep}). Since the YB operator $R$ is taken to
be an even Hecke symmetry of rank $p$, the image of the $p$-th
order antisymmetrizer
$$
{\cal A}^{(p)} := \rho_R(A^{(p)}):V^{\ot p}\to V^{\ot p}
$$
is one dimensional. So, the action of this projector on an
arbitrary basis element of $V^{\ot p}$ can be presented in the
form: \be {{\cal A}^{(p)}}e_{\langle 1|}\ot\dots \ot e_{\langle
p|}= e_{\langle 1|}\ot\dots \ot e_{\langle p|}\, v^{|1\dots
p\rangle}u_{\langle 1\dots p|},\label{u-v} \ee where components of
the tensors $u$ and $v$ belong to the field $\Bbb K $. Due to
(\ref{Y-ortog}) the normalization of $u$ and $v$ is fixed to be
$u_{\langle 1\dots p|}v^{|1\dots p\rangle} =1$. The images of
antisymmetrizers $A^{(m)}$ with $m>p$ are trivial: \be {\rm
Im}\,\rho_R(A^{(m)}) = 0\quad \forall m\ge p+1.\label{triv-im} \ee

In accordance with general theory each primitive idempotent
$Y^\lambda_{ii}$ generates a left or right ideal of the Hecke
algebra $H_m$ by means of left or right multiplications on all
possible elements of $H_m$. Equivalently, one can consider the
left or right {\em regular module} over $H_m$ which  as a vector
space is $H_m$ itself and left or right action consists in the
left or right multiplication by Hecke algebra elements. In these
terms the left or right ideals generated by primitive idempotents
will be irreducible left or right $H_m$ submodules. Let us denote
left and right submodules generated by $Y^\lambda_{ii}$ as
$M_\lambda^l(i)$ and $M_\lambda^r(i)$ respectively. As a direct
consequence of (\ref{id-dec}) one can write the following
decomposition: \be H_m = \bigoplus_{\lambda\vdash
m}\bigoplus_{i=1}^{d_\lambda}M_\lambda^r(i), \label{Dec} \ee and
a similar decomposition holds for left submodules
$M^l_\lambda(i)$.

Consider in more detail the structure of the component of
(\ref{Dec}) corresponding to a fixed partition $\lambda \vdash m$.
Let us denote \be \bar M_\lambda\equiv
\bigoplus_{i=1}^{d_\lambda}M_\lambda^r(i) \equiv
\bigoplus_{i=1}^{d_\lambda}M_\lambda^l(i). \label{M-bar} \ee The
set $\bar M_\lambda$ is a two sided submodule in $H_m$. Let us
dwell upon the question of linear basis in $\bar M_\lambda$. As
is known from the theory of Hecke algebra, $\bar M_\lambda$ is
isomorphic to the algebra $Mat_{d_{\lambda}}(\Bbb K)$ of
$d_{\lambda}\times d_{\lambda}$ matrices. Therefore in $\bar
M_\lambda$ one can find the set of $d_\lambda^2$ linear
independent quantities $Y^\lambda_{ij}$, which are some
polynomials in generators $\sigma_i$ of $H_m$ obtained from the
primitive idempotents $Y^\lambda_{ii}$. These quantities form a
linear basis in $\bar M_\lambda$ and obey the algebra of matrix
units\footnote{Matrix unit $E_{ij}$ $1\le i,j\le m$ of the matrix
algebra $Mat_m(\Bbb K)$ is the $m\times m$ matrix  with the only
nonzero entry at the intersection of the $i$-th row and the $j$-th
column which is equal to 1.} $E_{ij}$ which constitute the linear
basis in $Mat_{d_\lambda}(\Bbb K)$: \be
Y^\lambda_{ij}\,Y^\lambda_{kl} = \delta_{jk}\, Y^\lambda_{il},
\quad i,j,k,l = 1\dots d_{\lambda}. \label{Y-alg} \ee Note, that
for each fixed $i$ the elements $Y^\lambda_{ij}$ $1\le j\le
d_\lambda$ form a linear basis in the irreducible submodule
$M_\lambda^r(i)\subset \bar M_\lambda$ generated by the primitive
idempotent $Y^\lambda_{ii}$. And on the other hand, for each fixed
$j$ the elements $Y^\lambda_{ij}$ $1\le i\le d_\lambda$ form a
linear basis in $M_\lambda^l(j)$ generated by $Y^\lambda_{jj}$.
Therefore we have $\dim M^r_\lambda(i) = \dim M^l_\lambda(i)=
d_\lambda$.

Otherwise stated, if we arrange the elements $Y^\lambda_{ij}$ into
a rectangular $d_\lambda\times d_\lambda$ matrix, then its rows
will represent the linear basises of $d_\lambda$ right submodules
$M^r_\lambda(i)$, while the columns of the matrix will represent
the linear basises of left submodules $M^l_\lambda(i)$. The
diagonal entries of the matrix are the primitive idempotents
$Y^\lambda_{ii}$.

An important property of $Y^\lambda_{ii}$ consists in the
following fact. For any two idempotents $Y^\lambda_{ii}$ and
$Y^\lambda_{jj}$ corresponding to the same partition $\lambda$
there exists an invertible element of $H_m$ which transforms one
of these idempotents into the other one. This means that the
submodules $M_\lambda^r(i)$ (or $M_\lambda^l(i)$) with different
$i$ are isomorphic: they can be transformed into each other by
{\it left} (respectively {\it right}) multiplication on some
invertible element of Hecke algebra $H_m$. Note, that the whole
$\bar M_\lambda$, being a two-sided submodule in $H_m$, is
invariant with respect to such a multiplication.

In what follows we will not distinguish the isomorphic submodules
corresponding to a partition $\lambda$ and will consider them up
to an isomorphism. Then each submodule $M^r_\lambda(i)\subset \bar
M_\lambda$ generated by $Y^\lambda_{ii}$ can be treated as the
image of an irreducible $H_m$ module $M^r_\lambda$ (labelled only
by the partition $\lambda\vdash m$) w.r.t. the following $H_m$
module monomorphism: \be M_\lambda^r \rightarrow
M^r_\lambda(i)\subset \bar M_\lambda\subset H_m,\quad
\lambda\vdash m. \label{H-module} \ee In this sense the submodule
$\bar M_\lambda$ is isomorphic (as the $H_m$ module) to a direct
sum of $d_\lambda$ copies of the module $M^r_\lambda$:
$$
\bar M_\lambda\equiv
\bigoplus_{i=1}^{d_\lambda}M_\lambda^r(i)\cong
(M^r_\lambda)^{\oplus d_\lambda}.
$$
This formula reflects the well known fact that the regular
representation of a finite dimensional semisimple algebra
decomposes into a direct sum of irreducible modules and the
multiplicity of each module is equal to its dimension.

{\remark{\label{arb-emb} It is worth mentioning that the module
$M_\lambda^r$ can be mapped into $\bar M_\lambda$ in many
different ways and submodules $M_\lambda^r(i)$ (generated by a
fixed choice of primitive idempotents $\{Y^\lambda_{ii}\}$) are
only particular cases of all possible monomorphisms. Indeed, let
us take an arbitrary set of orthonormal $d_\lambda \times
d_\lambda$ projectors $P^{(i)}$ of the rank 1:
$$
P^{(i)}\cdot P^{(j)} = \delta^{ij}\,P^{(i)},\quad
\sum_{i=1}^{d_\lambda}P^{(i)} = {\id}.
$$
Then, as is evident from above relations and (\ref{Y-alg}), the
quantities $X^\lambda_{ii} = \sum_{r,s}P^{(i)}_{rs}
Y^\lambda_{rs}$ are also a set of primitive idempotents in $H_m $
leading to another decomposition of $\bar M_\lambda$ into a direct
sum of right submodules which will represent other possible
monomorphisms of $M_\lambda^r$ into $\bar M_\lambda$. }}

Now with each right\footnote{\label{right-sub}The choice of {\it
right} submodules is made since we prefer to use the {\it left}
action of $T\in {\rm End}(V)$ on $V$. Since an arbitrary element
of $M^r_\lambda(i)$ has the form $Y^\lambda_{ii}
f(\sigma_1\dots\sigma_m)$ ($f$ being a polynomial in the
generators of $H_m$), then with such a choice the projector
$\rho_R(Y^\lambda_{ii})$ will be the last element acting on
$V^{\ot m}$ in  formula (\ref{space}).} submodule $M^r_\lambda(i)
\subset H_m$ we can associate a space $V_\lambda(i)\subset V^{\ot
m}$ in the following way: \be V_\lambda(i) = {\rm
Im}\,\rho_R(M_\lambda^r(i)).\label{space} \ee Besides, we will
deal with a space $\bar V_\lambda$ which is the image of $\bar
M_\lambda$ (\ref{M-bar}): \be \bar V_\lambda = {\rm
Im}\,\rho_R(\bar M_\lambda)\subset V^{\ot m}.\label{bar-space} \ee

The spaces $V_\lambda(i)$ with different $i$ (and all other
$V_\lambda(X)\in \bar V_\lambda$ which are $\rho_R$-images of
other possible monomorphisms\footnote{See Remark~\ref{arb-emb}} of
$M_\lambda^r$ into $\bar M_\lambda$) are isomorphic as vector
spaces and we will not distinguish them. Instead, we will deal
with a space $V_\lambda$  which (like $M_\lambda^r$) gives rise to
a class of isomorphic embeddings $V_\lambda\hookrightarrow \bar
V_\lambda\subset V^{\ot m}$ and each $V_\lambda(i)$ (or any other
$V_\lambda(X)$ as well) is just a particular representative of
this class of isomorphic spaces.

{\remark{\label{zero-part} Formula (\ref{space}) has a nontrivial
meaning for $m\ge 1$, whereas at $m=1$ the only space $V_\lambda$ is
the space $V$ itself. For the future convenience we extend the
formula to the case $k=0$. Namely, we will take by definition
$V_{\lambda\vdash 0}=V_0\equiv \Bbb K$.}}

Now we can define the class of objects of the SW(V) category
generated by a finite dimensional space $V$ over a field $\Bbb K$
equipped with a Hecke symmetry $R$ (\ref{rmat})--(\ref{hec}).

{\definition{ To each fixed nonnegative integer $k\in\Bbb Z_+$
(for $k=0$ see remark \ref{zero-part} above) and each possible
partition $\lambda\vdash k$ we put into correspondence a space
$V_\lambda$  isomorphic to any $V_\lambda(i)$ in (\ref{space}).
The spaces $V_\lambda$, $\lambda\vdash k$, $k\in \Bbb Z_+$ are the
basic objects of the category SW(V). The whole class of objects of
the category is formed by direct sums of a finite number of basic
objects $V_\lambda$. Thus, the spaces $V_\la(i)$ can be treated as
the space $V_\la$ equipped with an embedding
$\{V_\lambda\hookrightarrow \bar V_\lambda\subset V^{\ot k}\}$.
\label{bas-obj}}}

{\remark{The objects $V_\la$ should be simple objects of our
category. However,  some objects $V_\la$ of the category  will be
identified with each other (in particular, $V_{(1^p)}$ and
$V_0$). So, finally simple objects of our category will be $V_\la$
modulo the mentioned identification. Up to this identification we
will sometimes use the notation $[V_\lambda]_k$ for the family of
embeddings $\{V_\lambda\hookrightarrow \bar V_\lambda\subset
V^{\ot k}\}$.}}

Let us turn now to the definition of morphisms of our category.

\subsection{\it Morphisms of the first kind and the structure of
quasitensor category\label{quasitensor}}

Let us denote $\Mor(U, V)$ the space of categorical morphisms
$U\to V$.

{\definition{\label{first-mor} The morphisms of the first kind are
defined as follows :\par i) The set $\Mor(V_\lambda, V_\lambda)$
for any basic object $V_\lambda$ contains only multiples of the
identical morphism. That is for any morphism $f:\,
V_\lambda\rightarrow V_\lambda$ we have by definition $f=a\,\id$
for some $a\in \Bbb K$.
\par ii)For an object embedded in $V^{\ot k}$
the morphisms of the first kind are represented by a set of
linear mappings:\par \be \forall\,k\in \Bbb N,\; \forall\,\tau\in
H_k,\quad\phi^k(\tau) = \rho_R(\tau): \quad V^{\ot k}\rightarrow
V^{\ot k}.\label{1-mor} \ee }}

Thus, being restricted to the set of subspaces $[V_\lambda]_k$
such a mapping $\phi^k$ sends each $V_\lambda(i)\in
[V_\lambda]_k$ to an isomorphic space or to zero
space\footnote{For example, $\rho_R(\sigma_i + q^{-1} \,{\rm
id}_{\! H})$ $1\le i\le k-1$ sends to zero the image of the
antisymmetrizer $A^{(k)}$.}. This means that nontrivial mappings
$\phi^k$ at most change the embedding of $V_\lambda$ into $V^{\ot
k}$ and therefore are multiples of identical morphism for the
basic objects. Among the morphisms $\phi^k$ there is no one which
would send $V_\lambda$ to $V_\mu$ with $\lambda\not=\mu$. A
different kind of morphisms is considered in subsection
\ref{red-proc}.

Our next step consists  in constructing a monoidal structure,
which allows us to "multiply" the objects of the category. This
means that we want to  define a covariant functor
$$
\otimes:\;{\rm Ob(SW)}\times{\rm Ob(SW)}\rightarrow {\rm Ob(SW)}
$$
with some associativity morphisms (for detail see \cite{CP,MacL}).
As such a functor we take the usual tensor product of linear
spaces with associativity morphisms to be identical. So, we have
to prove that the tensor product of any two objects of our
category is also an object, that is it can be decomposed into a
direct sum of basic objects. Evidently, one only needs to verify
this property for the tensor product of basic objects $V_\lambda$.

{\proposition{\label{lem-monoid}For given $\lambda\vdash n$ and
$\mu\vdash m$ the tensor product of two basic objects $V_\lambda$
and $V_\mu$ can be expanded into a direct sum of basic objects
$V_\nu$, $\nu\vdash (n+m)$:
$$
V_\lambda\otimes V_\mu = c_{\lambda\mu}^{\;\;\nu}\,V_\nu
$$
the coefficients $c_{\lambda\mu}^{\;\;\nu}$ being the
Littlewood-Richardson ones entering the formula for product of
Schur symmetric functions  $s_\lambda$ in $p$
variables\footnote{The detailed description of Schur functions and
related topics can be found in \cite{M}.}. }}

\smallskip

{\bf  Proof.}\ \  As is clear from the definition of basic objects
$V_\lambda$ of our category the structure of their tensor product
is controlled by that of modules $M_\lambda^r$ (\ref{H-module})
and the assertion of the proposition can be reformulated in terms
of these modules.

To do so let us consider an embedding of $V_\lambda\otimes V_\mu$
into $V^{\ot (n+m)}$ of the form
$$
V_\lambda\otimes V_\mu \rightarrow V_\lambda(i)\otimes
V_\mu(j)\subset V^{\ot (n+m)}
$$
for some fixed $i$ and $j$. As was defined in (\ref{space}) above,
the space $V_\lambda(i)$ is the image of the right $H_n$ submodule
$M_\lambda^r(i)$ under representation $\rho_R$ (\ref{Rep}).

Consider an embedding  ${\cal J}:\; H_n\times H_m \hookrightarrow
H_{n+m}$ which on generators looks as follows: \be {\forall
\sigma_i'\in H_n:\;\; {\cal J}(\sigma_i') = \sigma_i\in
H_{n+m}\quad 1\le i\le n-1 \atop \forall \sigma_j''\in H_m:\;\;
{\cal J}(\sigma_j'') = \sigma_{n+j}\in H_{n+m}\quad 1\le j\le
m-1}\label{Emb} \ee By construction ${\cal J}(H_n)$ and ${\cal
J}(H_m)$ form two mutually commuting Hecke subalgebras in
$H_{m+n}$.

Let $\lambda\vdash n$, $\mu\vdash m$ and $M_\lambda'(i)$,
$M_\mu''(j)$ be two right submodules in $H_n$ and $H_m$ generated
by the primitive idempotents $Y^\lambda_{ii}(\sigma')$ and
$Y^\mu_{jj}(\sigma'')$ correspondingly. The image of the tensor
product $M_\lambda'(i)\ot M_\mu''(j)$ under (\ref{Emb}) is a right
$H_n\times H_m$ module $M_\lambda^r(i)\ot M_\mu^r(j) \subset
H_{n+m}$ and obviously \be V_\lambda(i) \ot V_\mu(j) = {\rm
Im}\,\rho_R(M_\lambda^r(i)\ot M_\mu^r(j)).\label{V-rhoM} \ee

The right $H_{n+m}$ module induced from $M_\lambda^r(i)\ot
M^r_\mu(j)$ is reducible and can be decomposed into a direct sum
of irreducible right $H_{n+m}$ submodules $M^r_\nu(k)$. This
decomposition is the image (under monomorphism
$M_\lambda^r\rightarrow M_\lambda^r(i)$) of the following
relation between irreducible modules $M_\lambda^r$ and $M_\mu^r$:
\be M_\lambda^r\otimes M_\mu^r = c_{\lambda\mu}^{\;\;\nu}
\,M_\nu^r\label{prod-M} \ee where the coefficients
$c_{\lambda\mu}^{\;\;\nu}$ are equal to the multiplicity of
irreducible $H_{n+m}$ characters $\chi^\nu$ in the character
induced from $\chi^\lambda\times \chi^\mu$. Due to (\ref{V-rhoM})
the coefficients $c_{\lambda\mu}^{\;\;\nu}$ will also determine
the expansion of tensor product of two basic objects $V_\lambda$
and
 $V_\mu$.

Now we use the fact that at generic $q$ the Hecke algebra $H_m$ is
isomorphic to the group algebra $\Bbb K(S_m)$ of the $m$-th order
permutation group $S_m$ for all $m\ge 1$. Therefore, the
coefficients $c_{\lambda\mu}^{\;\;\nu}$ defining the multiplicity
of irreducible module $M_\nu^r$ in the tensor product
$M_\lambda^r\ot M_\mu^r$ are the same for the Hecke algebra and
for the algebra $\Bbb K(S_m)$. To complete the proof, note that
as is well known from the representation theory of symmetric
group (see for example \cite{M}) the corresponding multiplicities
coincide with the Littlewood-Richardson coefficients in the
product of symmetric Schur functions $s_\lambda$.\hfill
\rule{6.5pt}{6.5pt}

{\remark{The explicit calculation of the decomposition of
$M^r_\lambda(i)\otimes M^r_\mu(j)$ into a direct sum of $H_{n+m}$
submodules gives
$$
M_\lambda^r(i)\otimes M_\mu^r(j) =
\alpha_{\lambda\mu}^{\;\;\nu}(k(i,j))\,M_\nu^r(k(i,j)), \quad
\alpha_{\lambda\mu}^{\;\;\nu}(k(i,j))\in \Bbb Z_+.
$$
But the isomorphic submodules $M^r_\nu(k)$ with different $k$ are
images of the same module $M_\nu^r$ equipped with different
monomorphisms (depending on $i$ and $j$) $M_\nu^r\rightarrow \bar
M_\nu\subset H_{n+m}$. The sum $\sum_k
\alpha_{\lambda\mu}^{\;\;\nu}(k)$ does not depend on $i$ and $j$
and is equal to $c_{\lambda\mu}^{\;\;\nu}$ which defines the
structure of the tensor product (\ref{prod-M}). }}

Now we should convert the monoidal category SW(V) into a
quasitensor one (see \cite{CP}). This means, that we need to
define a set of {\em natural commutativity (or braiding)
isomorphisms} (\ref{rm}) for any pair $U, V\in {\rm Ob(SW(V))}$.
To be compatible with the monoidal structure the braiding
isomorphisms must satisfy the following property for any triple
$U,V, W\in {\rm Ob(SW(V))}$: \be R_{U, V\ot W} = ({\rm id}_V\ot
R_{U,W}) \circ (R_{U, V}\ot {\rm id}_W)\quad R_{V\ot W, U} =
(R_{V,U}\ot {\rm id}_W)\circ  ({\rm id}_V\ot
R_{W,U}).\label{step-prop} \ee Let us recall that the
associativity isomorphisms are taken to be identical.

Besides, considering two possible ways of transformation of $U\ot
V\ot W$ into $W\ot V\ot U$ one comes to the following condition:
\be (R_{V,W}\ot {\rm id}_U)\circ ({\rm id}_V \ot R_{U,W})\circ
(R_{U,V}\ot {\rm id}_W) = ({\rm id}_W\ot R_{U,V} )\circ
(R_{U,W}\ot {\rm id}_V)\circ ({\rm id}_U\ot R_{V,W}
).\label{braid} \ee

Let us begin with the simplest object --- the space $V$. As a
braiding isomorphism for $V\ot V$ we take the Hecke symmetry $R$.
Requirement (\ref{braid}) on $V^{\ot 3}$ transforms into the
Yang-Baxter equation (\ref{YB}) on $R$ which is satisfied by the
definition of $R$. The decomposition property (\ref{step-prop})
will be used as {\em a definition} of braiding isomorphisms for an
arbitrary tensor power of the space $V$. That is we take: \be
R_{V, V^{\ot k}} = R_kR_{k-1}\dots R_1\quad  R_{V^{\ot k}, V} =
R_1R_2\dots R_k.\label{Vk-V} \ee In the above formula a shorthand
notation $R_i\equiv R_{ii+1}$ is used. Introducing one more
notation for the chain of $R$ matrices
$$
R_{i\rightarrow j}\equiv\left\{{R_iR_{i+1}\dots R_j\quad {\rm if\
}j=i+n\ge i\atop R_iR_{i-1}\dots R_j\quad {\rm if\ }  j=i-n\le
i}\right.
$$
we can write a compact form for the braiding isomorphism in the
general case:
\begin{eqnarray}
R_{V^{\ot n}, V^{\ot m}} &\hspace*{-2.3mm}=& \hspace*{-2.3mm}
R_{m\rightarrow (m+n-1)}R_{(m-1)\rightarrow (m+n-2)}\dots
R_{1\rightarrow n}
\nonumber \\
&\hspace*{-2.3mm}\equiv &\hspace*{-2.3mm} R_{m\rightarrow
1}R_{(m+1)\rightarrow 2}\dots R_{(n+m-1)\rightarrow
n}.\label{Vm-Vn}
\end{eqnarray}
The two parts of this formula correspond to  two possible ways of
passing from $V^{\ot n}\ot V^{\ot m}$ to $V^{\ot m}\ot V^{\ot n}$.
Note, that all mappings \r{Vm-Vn} are morphisms of the first kind.

As for the braiding isomorphisms for $\Bbb K\ot V$ we take it to
be the usual flip: \be \Bbb K\ot V  = V\ot\Bbb K =
V.\label{cat:id} \ee

By definition of the objects of the category SW(V) any tensor
power $V^{\ot m}$ can be decomposed into a direct sum of basic
objects $V_\lambda$, $\lambda\vdash m$. Therefore by making use of
formula (\ref{Vm-Vn}) one can define the braiding isomorphism for
the tensor product of two arbitrary basic objects $V_\lambda\ot
V_\mu$, (and, therefore, for the tensor product of any couple of
objects) if we manage to prove that isomorphism (\ref{Vm-Vn})
does not "destroy" the structure of embeddings of $V_\lambda\ot
V_\mu$ into $V^{\ot n}\ot V^{\ot m}$. If it is the case then we
can take the restriction of (\ref{Vm-Vn}) onto ${\bar
V}_\lambda\ot {\bar V}_\mu$ as the braiding morphism for the
tensor product of $V_\lambda$ and $V_\mu$.

That is we have to show the following. Consider two arbitrary
embeddings ${\cal J}_i: V_\lambda\hookrightarrow V^{\ot n}$ and
${\cal J}_j: V_\mu\hookrightarrow V^{\ot m}$. Let ${\cal
J}_i(V_\lambda) = V_\lambda(i)$ and ${\cal J}_j(V_\mu) =
V_\mu(j)$. Let us recall that $V_\lambda(i)$ is defined by
(\ref{space}). We want to show that under isomorphism
(\ref{Vm-Vn}) one gets:
$$
V_\lambda(i)\ot V_\mu(j)\longrightarrow
R_{\lambda\mu}(V_\lambda(i)\ot V_\mu(j))
$$
where $R_{\lambda\mu}$ is an invertible operator from ${\rm
End}({\bar V}_\lambda\ot {\bar V}_\mu)$ (the space $\bar
V_\lambda$ is defined in (\ref{bar-space})) which {\it does not\/}
depend on $i$ and $j$.

{\proposition{For two given partitions $\lambda\vdash n$,
$\mu\vdash m$ and two arbitrary integers $1\le i\le d_\lambda$ and
$1\le j\le d_\mu$ consider the corresponding spaces $V_\lambda(i)$
and $V_\mu(j)$ as defined in (\ref{space}). Then under isomorphism
(\ref{Vm-Vn}) one has \be V_\lambda(i)\ot V_\mu(j) \longrightarrow
R_{\lambda\mu}(V_\lambda(i)\ot V_\mu(j)) \label{first} \ee where
the operator $R_{\lambda\mu}$ is defined  on $\bar V_\lambda \ot
\bar V_\mu$ by means of the following formula: \be R_{\lambda\mu}
= \left(\rho_R({\Bbb Y}^\mu)\ot\rho_R({\Bbb Y}^\lambda)
\right)\cdot R_{V^{\ot n}, V^{\ot m}}%
. \label{second} \ee Here $R_{V^{\ot n}, V^{\ot m}}$ is defined
by (\ref{Vm-Vn}) and
$$
\Bbb Y^\lambda\equiv\sum_{i=1}^{d_\lambda}Y^\lambda_{ii}
$$
is a central idempotent of the corresponding Hecke algebra.} }
\smallskip

{\bf Proof.}\  \ The proof is based on the embedding of $H_n\times
H_m$ into $H_{n+m}$ introduced in Proposition~\ref{lem-monoid}. In
accordance with (\ref{Emb}) and (\ref{V-rhoM}) we can write
$$
V_\lambda(i)\ot V_\mu(j) = {\rm
Im}\bigl\{\rho_R(Y^\lambda_{ii})\ot \rho_R(Y^\mu_{jj})\bigr\}.
$$
Considered on the arbitrary basis element of the space $V^{\ot
n}\ot V^{\ot m}$ the right hand side of the above formula reads:
\be e_{\langle 1|}\ot \dots \ot e_{\langle n+m|} \rightarrow
e_{\langle 1|}\ot \dots \ot e_{\langle n+m|}\, \bigl\{{\cal
Y}^\lambda_{ii}\ot {\cal Y}^\mu_{jj}\bigr\}^{|1\dots
n+m\rangle}_{\;\;\; \langle 1\dots n+m|} \label{mat-f} \ee where
the matrices ${\cal Y}^\lambda_{ii} = \rho_R(Y^\lambda_{ii})$ and
${\cal Y}^\mu_{jj} = \rho_R(Y^\mu_{jj})$ are some polynomials in
$R_1\dots R_{n-1}$ and $R_{n+1}\dots R_{n+m-1}$ respectively:
$$
{\cal Y}^\lambda_{ii} = {\cal Y}^\lambda_{ii}(R_1,\dots
,R_{n-1}),\quad {\cal Y}^\mu_{jj} = {\cal
Y}^\mu_{jj}(R_{n+1},\dots ,R_{n+m-1}).
$$
Applying the braiding isomorphism (\ref{Vm-Vn}) to $V^{\ot n}\ot
V^{\ot m}$ leads to the following result for formula (\ref{mat-f})
:
$$
e_{\langle 1|}\ot \dots \ot e_{\langle n+m|} \rightarrow
e_{\langle 1|}\ot \dots \ot e_{\langle n+m|} \, \bigl\{R_{V^{\ot
n},V^{\ot m}}\cdot({\cal Y}^\lambda_{ii}\ot {\cal
Y}^\mu_{jj})\bigr\}
$$
where the symbol $\,\cdot\,$ stands for the matrix multiplication.

Now one should take into account the following relations which are
direct consequence of (\ref{Vm-Vn}) and Yang-Baxter equation
(\ref{YB}):
$$
\begin{array}{ll}
R_{V^{\ot n},V^{\ot m}}\cdot R_i = R_{i+m}\cdot R_{V^{\ot
n},V^{\ot m}}&
1\le i \le n-1\\
\rule{0pt}{6mm} R_{V^{\ot n},V^{\ot m}}\cdot R_j = R_{j-n}\cdot
R_{V^{\ot n},V^{\ot m}}& n+1\le j \le n+m-1.
\end{array}
$$
Using these relations one gets:
\begin{eqnarray*}
R_{V^{\ot n},V^{\ot m}}\cdot {\cal Y}^\lambda_{ii}(R_1\dots
R_{n-1})&\hspace*{-2.3mm}
\ot& \hspace*{-2.3mm} {\cal Y}^\mu_{jj}(R_{n+1}\dots R_{n+m-1}) =\\
&&{\cal Y}^\mu_{jj}(R_{1}\dots R_{m-1})\ot {\cal
Y}^\lambda_{ii}(R_{m+1}\dots R_{m+n-1})\cdot R_{V^{\ot n},V^{\ot
m}}.
\end{eqnarray*}
This formula proves (\ref{first}). In order to find the form of
$R_{\lambda\mu}$ we observe that $\Bbb Y^\lambda = \sum_i
Y^\lambda_{ii}$ are {\it central} elements of $H_n$
($\lambda\vdash n$) and besides $Y^\lambda_{ii}\equiv
Y^\lambda_{ii} \Bbb Y^\lambda$ $\forall i,j$. Therefore
$$
({\cal Y}^\mu_{jj}\ot {\cal Y}^\lambda_{ii})\cdot R_{V^{\ot
n},V^{\ot m}}\equiv {\cal Y}^\mu_{jj}\ot {\cal
Y}^\lambda_{ii}\cdot\bigl(\rho_R(\Bbb Y^\mu)\ot \rho_R(\Bbb
Y^\lambda)\cdot  R_{V^{\ot n},V^{\ot m}}\bigr).
$$
where $\rho_R(\Bbb Y^\lambda)$ is the projector onto $\bar
V_\lambda\subset V^{\ot (n+m)}$. Thus, we come to form
(\ref{second}) of $R_{\lambda\mu}$.\hfill  \rule{6.5pt}{6.5pt}

So, the operator $R_{\lambda\mu} = \rho_R(\Bbb Y^\mu)\ot
\rho_R(\Bbb Y^\lambda)R_{V^{\ot n},V^{\ot m}}$ does not depend on
a concrete embedding of $V_\lambda \ot V_\mu \in {\rm Ob(SW(V))}$
into $V^{n}\ot V^m$ (i.e., it does not depend on indices $i,j$)
and represents a braiding isomorphism for $V_\lambda\ot V_\mu$.

\subsection{\it Morphisms of the second kind and reduction
procedure\label{red-proc}}

As was mentioned above (see (\ref{triv-im})) the image of the
antisymmetrizer $A^m$ with $m > p$ is identical zero in $V^{\ot
m}$. It can be shown that the same is true for any ${\rm
Im}\,\rho_R(M^r_\lambda)$ in case if the height
$\ell(\lambda)>p$. This means that the basic objects of our
category $V_\lambda$  are labelled by partitions with restricted
height: $V_\lambda \not\equiv 0\Leftrightarrow \ell(\lambda) \le
p$.

However, there is  another consequence of (\ref{u-v}) --
(\ref{triv-im}) which concerns the objects labelled by partitions
with the height equal to $p$. Consider a partition $\lambda\vdash
(p+m)$ for some nonnegative integer $m$ such that $\ell(\lambda) =
p$. Let $\mu$ denote a partition of $m$ which is obtained from
$\lambda$ by striking the first column out of the Young diagram
corresponding to $\lambda$. The parts of $\lambda$ and $\mu$ are
connected by the relation $\mu_i = \lambda_i - 1$, $\forall
\lambda_i\not= 0$. As immediately follows from
Proposition~\ref{lem-monoid} and from the above remark about the
maximal height of $\lambda$ the following decomposition takes
place: \be V_\lambda\cong V_{(1^p)}\ot V_\mu\label{eq:sp-dec} \ee
where $V_{(1^p)}$ stands for the one dimensional space labelled
by one-column diagram with $p$ boxes: $\lambda = (1^p)$.

Isomorphisms (\ref{eq:sp-dec}) allow us to define the so called
{\it reduction procedure} and to introduce morphisms of the second
kind in our category, namely, those identifying $V_\lambda$ and
$V_\mu$. For this purpose consider a mapping $\psi$ sending the
one-dimensional space $V_{(1^p)}$ into the field $\Bbb K$.
Evidently, one only needs to fix the action of the mapping $\psi $
on a basis vector of $V_{(1^p)}\cong {\rm Im}\,\rho_R(A^{(p)})$.
In the general form such a mapping looks like \be \psi(e_{\langle
1|}\ot \dots \ot e_{\langle p|}\,v^{|1\dots p\rangle}) = 1\in\Bbb
K.\label{psi} \ee The choice of $1$ in the above formula does not
restrict generality of our consideration.

Let us suppose that $\psi\in\Mor(V_{(1^p)}, \K)$. Then taking
into account \r{nat} and the fact that the identity operator is a
categorical morphism
 we conclude that the following diagrams must be
commutative:
$$
\begin{array}{ccc}
{V_{(1^p)}\ot V_\mu}&\stackrel{\psi\ot {\rm
id}}{\longrightarrow}&{\Bbb
K\ot V_\mu}\\
\vcenter{\llap{$\scriptstyle R_{V_{(1^p)},V_\mu}$}}\Big\downarrow
&
&\Big\|\\
V_\mu\ot V_{(1^p)} & \stackrel{{\rm id}\ot\psi
}{\longrightarrow}&V_\mu\ot \Bbb K
\end{array}
\qquad\qquad
\begin{array}{ccc}
V_\mu\ot V_{(1^p)} & \stackrel{{\rm id}\ot\psi
}{\longrightarrow}&V_\mu\ot
\Bbb K\\
\vcenter{\llap{$\scriptstyle R_{V_\mu, V_{(1^p)}}
$}}\Big\downarrow &
&\Big\|\\
{V_{(1^p)}\ot V_\mu}&\stackrel{\psi\ot {\rm
id}}{\longrightarrow}&{\Bbb K\ot V_\mu}
\end{array}
$$

Consider a particular case of these diagrams by putting $V_\mu=V$
(if they are commutative in this case the same will be true for
any $\mu$). Then by passing to the basis we have \be V_{(1^p)}\ot
V:\quad e_{\langle 1|}\ot\dots\ot e_{\langle p|}\ot e_{\langle
p+1|}\,v^{|1\dots p\rangle} \stackrel{\psi\ot {\rm
id}}{\longrightarrow} 1\ot e_{\langle p+1|} = e_{\langle p+1|}\ot
1. \label{one-way} \ee

According to another way in the left diagram we should first apply
the braiding morphism $R_{V_{(1^p)},V}$. Thus, we have
\begin{eqnarray*}
R_{V_{(1^p)},V}:\;\;e_{\langle 1|}\ot\dots\ot e_{\langle p|}\ot
e_{\langle p+1|}\,
v^{|1\dots p\rangle}&&\longrightarrow\\
&&e_{\langle 1|}\ot\dots\ot e_{\langle p|}\ot e_{\langle
p+1|}\,R_1\dots R_p\,v^{|1\dots p\rangle}.
\end{eqnarray*}
In order to simplify this expression we need the following useful
relations (see \cite{5avt}):
\begin{eqnarray}
&&R_1\dots R_p\,{\cal A}^{(p)} = (-1)^{p-1}\,qp_q\,{\cal
A}^{(2,p+1)}
{\cal A}^{(p)} \nonumber\\
&&R_p\dots R_1\,{\cal A}^{(2,p+1)}= (-1)^{p-1}\,qp_q\,{\cal
A}^{(p)}
{\cal A}^{(2,p+1)}\label{A-rel}\\
&& {\cal A}^{(p)}{\cal A}^{(2,p+1)}{\cal A}^{(p)} = p_q^{-2} {\cal
A}^{(p)}\ot {\rm id}_{p+1}. \nonumber
\end{eqnarray}
In these formulas the antisymmetrizer ${\cal A}^{(2,p+1)}$ has the
same form as ${\cal A}^{(p)}$ but it depends on $R_2,\dots ,R_{p}$
instead of $R_1,\dots ,R_{p-1}$.

Now in virtue of (\ref{A-rel}) we obtain: \be R_1\dots
R_p\,v^{|1\dots p\rangle} = v^{|2\dots
p+1\rangle}\,N^{|1\rangle}_{\;\langle p+1|}\, \label{trans} \ee
 where the matrix
$N$ is defined as follows: \be N^i_{\;j}\equiv (-1)^{p-1} q p_q
\,u_{a_2\dots a_p j}v^{ia_2\dots a_p}  \quad{\rm or}\quad
N^{|1\rangle}_{\;\langle 1'|} = (-1)^{p-1} q p_q \,u_{\langle
2\dots p 1'|}v^{|12\dots p\rangle}.\label{def:M-N} \ee

Therefore the composition $({\rm id\ot \psi})\circ
R_{V_{(1^p)},V}$ gives the following result:
$$
e_{\langle 1|}\ot\dots\ot e_{\langle p|}\ot e_{\langle
p+1|}\,v^{|1\dots p\rangle} \rightarrow e_{\langle 1|}\ot\dots\ot
e_{\langle p|}\ot e_{\langle p+1|}\,v^{|2\dots p+1\rangle}\,
N^{|1\rangle}_{\;\langle p+1|}\rightarrow e_{\langle
1|}N^{|1\rangle}_{\;\langle p+1|}\ot 1
$$
which obviously differs from (\ref{one-way}). These results can be
made compatible iff $N$ is a scalar matrix (i.e., it is
proportional to the $n\times n$ unit matrix $I$). If it is so,
then by multiplying the braiding $R$ by an appropriate factor we
can achieve the commutativity of the left diagram above. Indeed,
if we replace the braiding $R$ by $aR$, $a\in\Bbb{K}$, the factor
$a^p$ appears in the right hand side of \r{trans}. Choosing
properly the factor $a$ we can obtain $a^p N=I$.

Applying the same procedure for the space $V\ot V_{(1^p)}$ (the
right diagram above) we find:
$$
e_{\langle 1|}\ot e_{\langle 2|}\ot\dots\ot  e_{\langle
p+1|}\,v^{|2\dots p+1\rangle} \longrightarrow 1\ot e_{\langle
p+1|} M^{|p+1\rangle}_{\;\;\langle 1|}
$$
where\footnote{Note that the matrices $N$ and $M$ differ by a
factor from those considered in \cite{G}.}
$$
M^i_{\;j}\equiv (-1)^{p-1} q p_q\,u_{j a_2\dots a_p}v^{a_2\dots
a_p i} \quad{\rm or}\quad M^{|1'\rangle}_{\;\langle 1|}\equiv
(-1)^{p-1} q p_q\,u_{\langle 12\dots p|}v^{|2\dots p 1'\rangle}.
\eqno{(\ref{def:M-N}')}
$$

For the same reason we assume the matrix $M$ to be scalar. Then
the right diagram can also be made commutative by a proper
rescaling of $R$. However, if we want to obtain the unit matrix
instead of the matrices $N$ and $M$ simultaneously  we have to
impose one more condition: $M=N$.

Therefore from now on we will suppose the Hecke symmetry $R$ to
satisfy the relation \be
 M = N=a I, \quad a\in \Bbb K^\times .\label{restr}
\ee Thus, the mapping $\psi$ becomes a morphism of the category
SW(V) after a proper renormalization  of the braiding morphisms.

This normalization factor is easy to find. Indeed, using relations
(\ref{A-rel}) (valid for {\it any} Hecke symmetry of the rank
$p$) we find \be M\cdot N = q^2 I.\label{N-inv-M} \ee So, if
relation \r{restr} is satisfied  then by virtue of \r{N-inv-M} we
have $a=\pm q$ and hence \be N=M=\pm q I.\label{fix-MN} \ee

Therefore by assuming $\psi$ to be a morphism and the braiding to
be natural we should pass from the initial Hecke symmetry $R$ to
the braiding \be \bar R= (\pm q)^{-\frac{1}{p}}R.\label{renor}
\ee On the higher tensor powers this renormalization gives rise
to a renormalized braiding as well \be \bar R_{V^{\ot n}, V^{\ot
m}} = (\pm q)^{-\frac{mn}{p}}R_{V^{\ot n}, V^{\ot
m}}.\label{renorm} \ee The mapping $\psi$ (as well as its inverse
$\psi^{-1}$) will be called a morphism of the second kind.

Now let us compare our approach to introducing categorical
morphisms with that making use of the RTT algebra. Recall that
this algebra denoted $\cal T$ is generated by $n^2$ elements
$T^i_{\;j}$ subject to the relations \cite{FRT}: \be R_{12}T_1T_2
= T_1T_2R_{12}.\label{RTT} \ee It becomes a bialgebra being
equipped with a coproduct $\Delta$ and a counit $\varepsilon$ as
follows:
$$
\Delta(T^i_{\;j}) = T^i_{\;k}\otimes T^k_{\;j}\qquad
\varepsilon(T^i_{\;j}) = \delta^i_{\;j}.
$$
Define the right comodule structure $\delta_r:\;V\rightarrow
V\otimes \cal T$ on the space $V$ as follows:
$$
\delta_r(v) = e_k\otimes T^k_{\;\,i}v^i,\quad \forall
\;v=v^ie_i\in V.
$$
Such a coaction is extended to $V^{\ot k}$ in the obvious way: \be
\delta_r: V^{\ot k}\rightarrow V^{\ot k}\ot {\cal T}\quad
e_{\langle 1|}\ot\dots \ot e_{\langle k|}\rightarrow e_{\langle
1|}\ot\dots \ot e_{\langle k|}\ot T_1\dots T_k\label{def:comod}
\ee and all the comodule properties are easily verified. Then each
$V_\lambda(i)$ $\lambda\vdash k$ turns out to be an invariant
subcomodule in $V^{\ot k}$ since the coaction $\delta_r$ commute
with the action of $\rho_R(M_\lambda^r(i))$ on $V^{\ot k}$:
$$
\forall k,\;\forall\lambda\vdash k,\quad \rho_R(M_\lambda^r)\,
T_1T_1\dots T_k = T_1T_1\dots T_k\,\rho_R(M_\lambda^r).
$$
This equation is a direct consequence of defining relations
(\ref{RTT}) and explicit formula for pri\-mi\-tive idempotents of
$H_k$  which expresses each $Y^\lambda_{ii}$ as a polynomial in
$\sigma_1, \sigma_2,\dots, \sigma_{k-1}$.

Thus, any object of the category $\SW$ can be equipped with a
${\cal T}$-coaction. This structure is often used in order to
define categorical morphisms: one says that a map $U\to V$ is a
morphism if it commutes with this coaction. Let us analyze what it
entails being applied to $\psi$. For this purpose consider {\it a
quantum determinant} $\det_qT$ of bialgebra (\ref{RTT}) defined
as follows \cite{FRT}: \be {\det}_qT\equiv u_{\langle 1 2 \dots
p|}T_1T_2\dots T_p v^{|1 2 \dots p\rangle}.\label{q-det} \ee

The commutation relations of $\det_qT$ with generators of
(\ref{RTT}) read (cf. \cite{G}): \be {\det}_qT\cdot T =
(N^{-1}T\,N)\cdot {\det}_qT \quad {\Leftrightarrow} \quad
{\det}_qT\cdot (NT) = (TN)\cdot {\det}_qT \label{NTD} \ee or due
to (\ref{N-inv-M})
$$
{\det}_qT\cdot T = (MTM^{-1})\cdot {\det}_qT \quad
{\Leftrightarrow} \quad {\det}_qT\cdot (TM) = (MT)\cdot {\det}_qT.
$$

If we want the mapping $\psi$ (\ref{psi}) to commute with the
comodule structure we have to introduce one more condition: \be
{\det}_q T = 1. \label{det-cen} \ee But such a relation is
compatible with the algebraic structure only iff $\det_q T$ is a
central element of algebra (\ref{RTT}). By virtue of (\ref{NTD})
this means that the matrix $N$ (and hence $M$) must be scalar.

Formally, the condition "$N$ and $M$ are scalar" is weaker than
"$N$ and $M$ are scalar and equal to each other". It is not clear
whether there exist Hecke symmetries such that the matrices $N$
and $M$ are scalar but not equal to each other. Nevertheless, the
family of even Hecke symmetries satisfying (\ref{fix-MN}) is
sufficiently large. It will be shown by one of the authors (R.L)
in a separate paper.

Let us also remark that for the quasiclassical Hecke symmetries
coming from the universal $R$ matrix of QG $U_q(sl(n))$ condition
(\ref{restr}) holds true. The quantum determinant $\det_qT$ and
the unit element $1$ generate the center of (\ref{RTT})
\cite{FRT}. But in general (even in the quasiclassical case)
$\det_qT$ is not central.

\subsection{\it The structure of rigid category}

To convert a (quasi)tensor category ${\cal C}$ into a {\it rigid}
one we must specify the following data \cite{CP}:
\begin{enumerate}
\item[a.] A mapping $*:\;{\rm Ob}({\cal C})\rightarrow  {\rm Ob}({\cal C})$
which for any $U\in {\rm Ob}({\cal C})$ put into correspondence
its {left dual} $U^*$
\item[b.] For any pair of dual objects $U$ and $U^*$ there exist evaluation
$\ev{U}$ and coevaluation $\mpi{U}$ morphisms of ${\cal C}$
$$
\mpi{U}:\; \Bbb K \longrightarrow U\ot U^*,\qquad \ev{U}:\; U^*\ot
U\longrightarrow \Bbb K
$$
such that the following diagrams are commutative: \be
\begin{array}{ccc}
{U}&\stackrel{\mpi{U}\ot {\rm id}}{\longrightarrow}&{(U \ot U^*)\ot U}\\
\vcenter{\llap{{\rm id}}}\Big\downarrow & &\Big\downarrow
\vcenter{\rlap{\rm id}}\\
U& \stackrel{{\rm id}\ot \ev{U}}{\longleftarrow}& {U\ot (U^*\ot
U)}
\end{array}
\qquad\qquad
\begin{array}{ccc}
{U^*}&\stackrel{{\rm id}\ot \mpi{U}}{\longrightarrow}&{ U^*\ot
(U\ot
U^*)}\\
\vcenter{\llap{{\rm id}}}\Big\downarrow & &\Big\downarrow
\vcenter{\rlap{\rm id}}\\
U^*& \stackrel{\ev{U}\ot {\rm id}}{\longleftarrow}& {(U^*\ot U)
\ot U^*}
\end{array}
\label{comm-diag} \ee
\end{enumerate}

Let us begin with the case $U=V$. By our construction of the
morphism $\psi$ it is evident that the dual space $V^*$ is
nothing but $\Lambda^{p-1}(V)$. So, in order to satisfy
\r{comm-diag} we should only fix a convenient basis converting
$\Lambda^{p-1}(V)$ into a left dual space.

{\definition{\label{def:dual} The space $\Lambda^{p-1}(V)$
equipped with the basis \be e^i = e_{a_2}\ot \dots\ot
e_{a_{p}}\,v^{ia_2\dots a_{p}} \quad{\rm or}\quad e^{|1\rangle} =
e_{\langle 2|}\ot \dots\ot e_{\langle p|}\,v^{|12\dots
p\rangle}.\label{du-bas} \ee will be called a left dual space to
the space $V$. The morphisms $\ev{V}$ and $\mpi{V}$ are defined
as follows (on basis vectors): \be
\begin{array}{l@{\,}cl}
\ev{V}&:& e^i\ot e_j \rightarrow \delta^i_{\; i}\,1\\
\rule{0pt}{6mm} \mpi{V}&:& 1 \longrightarrow e_i\ot e^i.
\end{array}
\label{rig-mor} \ee }}

In order to justify this definition we have to show the following.
{\proposition{Mappings (\ref{rig-mor}) are morphisms of the SW(V)
category and they satisfy pro\-per\-ty (\ref{comm-diag}).}}

\smallskip

{\bf Proof.}\ \ The fact that $\mpi{V}$ is a morphism is evident
since
$$
\mpi{V} = \psi^{-1}.
$$
To show that $\ev{V}$ is a morphism let us consider the mapping
$$
\varphi=(-1)^{p-1}p_q \psi\circ\mathcal{A}^p:\quad V^{\ot p}\to
\Bbb{K}.
$$
It is a morphism by construction. In a basis form we get
$$
\varphi:\; e_{i_1}\ot...\ot e_{i_p}\to (-1)^{p-1}p_q
u_{i_1...i_p}.
$$

Applying this morphism to the element
$$
e^i\ot e_j = v^{i a_2...a_p}e_{a_2}\ot...\ot e_{a_p}\ot e_j
$$
we get
$$
\varphi(e^i\ot e_i)=(-1)^{p-1}p_q v^{i a_2...a_p}u_{a_2...a_pj}
=q^{-1}N_j^i= \pm\delta_j^i.
$$
Thus, the mapping $\ev{V}$ is nothing but the morphism
$\pm\varphi$ restricted on $\Lambda^{p-1}(V)\ot V$.

Diagrams \r{comm-diag} are obviously commutative.\hfill
\rule{6.5pt}{6.5pt}

Let us define now left duals to other simple objects of our
category. First, introduce a dual object to $V^{\ot m}$ by
putting $(V^{\ot m})^*=(V^*)^{\ot m}$.

{\proposition{The space $(V^*)^{\ot m}$ is dual to the space
$(V^{\ot m})$ being equipped with the mappings
$$
\ev{U}\,:{V^*}^{\ot m}\ot {V}^{\ot m}\longrightarrow {\Bbb
K},\qquad \mpi{U}\,: \Bbb K \longrightarrow V^{\ot m}\ot
{V^*}^{\ot m}
$$
which are defined to be: \be
\begin{array}{l@{\,}c@{\quad}l}
\ev{U}&:& e^{|m\rangle}\ot \dots \ot e^{|1\rangle}\ot e_{\langle
1|}\ot \dots e_{\langle m|} \rightarrow
\delta^{|1\rangle}_{\;\langle 1|}\dots \delta^{|m\rangle}_{\;\langle m|}\\
\rule{0pt}{6mm} \mpi{U}&:& 1_{\Bbb K}\rightarrow e_{\langle
1|}\ot\dots\ot e_{\langle m|} \ot e^{|m\rangle}\ot \dots \ot
e^{|1\rangle}
\end{array}.
\ee  }}

\smallskip

{\bf Proof}\ \  is obvious.

Now we are ready to define the dual $V^*_\lambda$ and
corresponding morphisms $\mpi{V_\lambda}$ and $\ev{V_\lambda}$ for
arbitrary basic object $V_\lambda$ of the  category  $\SW$.

Fix a basic object $V_\lambda$, $\lambda\vdash m$, and consider
its particular embedding into $V^{\ot m}$ in the form
(\ref{space}):
$$
V_\lambda(i) = {\rm Im}\,\rho_R(Y^\lambda_{ii}).
$$
This space is spanned by the following vectors
$$
e_{\langle 1|}\ot \dots \ot e_{\langle m|}\,\cdot({\cal
Y}^\lambda_{ii})^{|1\dots m\rangle}_{\;\;\langle1\dots m|}
$$
where ${\cal Y}^\lambda_{ii}$ is a matrix of $Y^\lambda_{ii}$ in
the representation $\rho_R$.

Let us define a space $V^*_\lambda(i)\subset {V^*}^{\ot m}$ as a
linear span of the form: \be V^*_\lambda(i) =
\{w^*:\;w^*=w_{\langle 1\dots m|}\,({\cal Y}^\lambda_{ii})\cdot
e^{|m\rangle}\ot\dots \ot e^{|1\rangle}\}\label{space*} \ee where
$w_{\langle 1\dots m|}$ is an arbitrary tensor with components
from $\Bbb K$. It is important, that the order of spaces in this
formula is reversed comparing with the preceding relation while
the matrix ${\cal Y}^\lambda_{ii}$ is the same.

Finally, we have the following.

{\proposition{The space $V^*_\lambda$  being equipped with the
morphisms $\ev{V_\lambda}$ and $\mpi{V_\lambda}$ of the form \be
\begin{array}{l@{\,}cl}
\ev{V_\lambda}&:& V_\lambda^*\ot V_\lambda \longrightarrow {\cal
Y}^\lambda_{ii}\\
\rule{0pt}{6mm} \mpi{V_\lambda}&:&  1_{\Bbb K}\longrightarrow
e_{\langle 1|}\ot\dots e_{\langle m|}\ot {\cal Y}^\lambda_{ii}
\cdot e^{|m\rangle}\ot \dots \ot e^{|1\rangle}.
\end{array}
\ee becomes left dual to the space $V_\lambda$ (provided that
$V_\lambda$ is realized as $V_\lambda(i)$).}}

Let us note that in fact the space $V_\lambda^*$ is nothing but
$$
V_\mu,\, \,\,\mu=(\lambda_1, \lambda_1-\lambda_{p-1},
\lambda_1-\lambda_{p-2},...,\lambda_1-\lambda_2)
$$
equipped with an embedding in the space ${V^*}^{\ot m}$. This
shows that we do not need to construct left dual to the space
$V_\lambda^*$ because it is just the object $V_\lambda$.

Nevertheless, the problem of an explicit pairing
$$
V_\lambda\ot V_\lambda^*\to \Bbb{K}
$$
converting $V_\lambda^*$ into the right dual to $V_\lambda$ (and
$V_\lambda$ in the left dual to $V_\lambda^*$) is very important.
In order to describe this pairing in a particular case
$V_\lambda=V$  let us calculate the braiding of $V$ and
$\Lambda^{p-1}(V)$ expressed via the basis $\{e^i\}$. Note, that
we use the renormalized braiding $\bar R$ (\ref{renor}).

{\proposition{\label{star-braid} We have \be \bar R_{V, V^*}:\;\;
e_i\ot e^j \longrightarrow e^r\ot e_s\,\bar
Q^{js}_{\;ir}\label{V-V*} \ee where in the given basis $\{e_i\}$
the operator $\bar Q$ satisfies the relation: \be \bar R^{i
a}_{\;j b}\, \bar Q^{b k}_{\; a l} =
\delta^i_{\;l}\delta^k_{\;j}\quad\Leftrightarrow \quad \bar Q^{i
a}_{\;j b}\,\bar R^{b k}_{\; a l} = \delta^i_{\;l}\delta^k_{\;j}
\label{R-Q}. \ee }}

Note, that if such $\bar Q$ exists then the corresponding
braiding $\bar R$ is usually called {\it invertible by  column}.
\smallskip

{\bf Proof.}\ \  According to definition (\ref{du-bas}) for $e^i$
and (\ref{renorm}) for the renormalized braiding we get (omitting
the obvious signs of tensor product):
$$
e_i\ot e^j = e_i\,e_{a_2}\dots e_{a_p}\,v^{ja_2\dots a_p}
\stackrel{\mbox{\footnotesize $\bar R_{V,V^{\ot
(p-1)}}$}}{\longrightarrow} e_{s_1}\dots e_{s_p}(\pm
q)^{\frac{1-p}{p}}(R_{p-1}\dots R_1)^{s_1\dots s_p}_{\;\;
ia_2\dots a_p}v^{ja_2\dots a_p}.
$$

Now from (\ref{A-rel}) it is easy to find
$$
(\pm q)^{\frac{1-p}{p}}(R_{p-1}\dots R_1)^{s_1\dots s_p}_{\;\;
ia_2\dots a_p}v^{ja_2\dots a_p} = \bar
Q^{js_p}_{\;ir}\,v^{rs_1\dots s_{p-1}}.
$$
Substituting this into the previous relation and applying again
the definition (\ref{du-bas}) of the dual basis we find result
(\ref{V-V*}).\hfill \rule{6.5pt}{6.5pt}

Now we introduce a pairing between $V$ and $V^*$ by putting
$$
\bar\ev{V}=\ev{V}\circ R_{V, V^*}:\quad V\ot V^*\to \Bbb{K}.
$$
It is this operator which plays the central role in the next
section where a categorical trace in the space $\End(V)$ will be
introduced.

Completing this Section let us discuss a meaning of the property
that the morphism $\psi$ is coordinated with the braidings in the
sense of (\ref{nat}). Let $V$ be the basic object of the category
SW(V). Assuming the category to be rigid we can identify
$\End(V)\cong V\ot \Lambda^{(p-1)}(V)$ (we consider the left
morphisms space). Then the usual operator product $\mu:
\End(V)^{\ot 2}\to \End(V)$ is nothing but the pairing
$\Lambda^{(p-1)}(V)\ot V\to \K$. Moreover, the property
$$ R\,\mu_{12}=\mu_{23}\,R_{12}\,R_{23},\,\,\,R\,\mu_{23}=
\mu_{12}\,R_{23}\,R_{12},\,\,\, R=R_{\End(V), \End(V)}$$ which
appears in numerous applications is satisfied.  Let us note that
these relations are not invariant w.r.t. a rescaling $R\to
a\,R,\,\,a\not=0$ while the condition (\ref{det-cen}) is.

Also remark that often a problem arises to check that a given map
$\rho: U\to V,\,\,U,V\in \Ob(\B)$ is a categorical morphism. Then
by putting in (\ref{nat}) $f=\rho,\, g=\id$ or $f=\rho,\, g=\rho$
we get necessary conditions very useful in practice.

\section{Additive-multiplicative functional and dimension}

In this section we find all possible additive and multiplicative
(a-m) functionals on objects of Schur-Weyl category SW(V) and
suggest a definition of the trace in ${\rm End }(V)$ as morphism
(\ref{tr}) in the category such that the corresponding dimension
is an a-m functional. As was pointed out in the Introduction, the
requirement that the trace should be a morphism in general is not
a trivial one. In particular, this requirement gives rise to the
notion of {\it quantum trace} in the category of finite
dimensional modules over a special class of quasitriangular Hopf
algebras --- the so called ribbon Hopf algebras. In our approach
the morphism property of trace can be obtained without any
additional Hopf structure.

Let us recall that a mapping $f:\;{\rm Ob(SW(V))}\rightarrow \Bbb
K$ will be called {\em an ad\-di\-tive-mul\-ti\-pli\-ca\-tive
functional} (a-m functional) if it possesses the following
property: \be \forall\; U,V\in {\rm Ob(SW(V))}:\quad f(U\oplus V)
= f(U)+f(V),\quad f(U\otimes V) = f(U)f(V).\label{am-prop} \ee
{\remark{Introducing the Grothendieck semiring of the category
SW(V) and making use of the fact that our category is semisimple
(i.e. any its object can be decomposed into a direct sum of simple
ones) we conclude that each a-m functional on the category is in
fact defined on the Grothendieck semiring.}}

The fact that we are working with a Hecke symmetry of the rank $p$
leads to an important consequence for the set of a-m functionals.
Namely we will show that the value of any a-m functional $f$ on a
basic object $V_\lambda$ (and, hence, on an arbitrary object of
SW(V)) is completely determined by $(p-1)$
numbers\footnote{\label{p-dim} Due to (\ref{cat:id}) and reduction
morphism $\psi$ (\ref{psi}) we have to put for any a-m functional
$f(\Bbb K) = f(V_{(1^p)}) = 1$.} $f(V_{(1^k)})$ $k=1,\dots ,p-1$.

{\proposition{\label{AM-func} Given an a-m functional $f$,
construct a $p$-th order polynomial in a formal variable $t$ of
the form: \be \phi(t) = t^p + f_1t^{p-1}\dots +f_{p-1}t+1, \quad
f_k\equiv f(V_{(1^k)}), \;\; k=1,\dots,p-1. \label{polynom} \ee
Let numbers $(-\alpha_i)\in \Bbb K$ be the roots of $\phi(t)$,
that is $\phi(-\alpha_i) = 0$ $i = 1,\dots ,p$. Then for any basic
object $V_\lambda$ of SW(V) the following relation holds: \be
f(V_\lambda) = s_\lambda(\alpha_1, \alpha _2,\dots,
\alpha_p)\label{val} \ee where $s_\lambda(x_1,\dots,x_p)$ is the
symmetric Schur function in $p$ variables\footnote{Let us recall,
that due to (\ref{cat:id}) all the partitions $\lambda$ labelling
the objects of our category has the height not greater than $p$
therefore the right hand side of equation (\ref{val}) is correctly
defined for any object of SW(V).} corresponding to the partition
$\lambda$. }}
\smallskip

{\bf Proof.}\ \ First of all let us prove that the numbers $f_k =
f(V_{(1^k)})$ do satisfy (\ref{val}). In accordance with the
definition of $\alpha_i$ we can write
$$
\phi(t) = (t+\alpha_1)(t+\alpha_2)\dots (t+\alpha_p).
$$
Therefore, as immediately follows from this relation, the
coefficient $f_k$ is the $k$-th elementary symmetric function
$e_k$ in variables $\alpha_i$ \cite{M}:
$$
f_k = \sum_{i_1<i_2<\dots<i_k}\alpha_{i_1}\alpha_{i_2}\dots
\alpha_{i_k} = e_k(\alpha_1,\dots ,\alpha_k).
$$
But since $s_{(1^k)} = e_k$ we conclude that the assertion of the
proposition is valid for $f_k$.

The fact that the quantities $s_\lambda(\alpha_1, \alpha _2,\dots,
\alpha_p)$ can be considered as values of an a-m functional
follows immediately from the properties of Schur functions and
from Proposition~\ref{lem-monoid}. In order to prove \r{val} we
should only check that once the quantities $f(V_\lambda)$,
$\lambda=(1^k)$, are given then all quantities $f(V_\mu)$ for all
other partitions are uniquely defined. This can be shown by
induction in couples $(m,k)$ where $m$ is the number of columns
in $\lambda$ and $k$ is the number of boxes in the last column.
Namely, using the implication
$$
f(V_\lambda\ot V_\mu) =f(V_\lambda)f(V_\mu) =
c_{\lambda\mu}^{\;\;\nu}\,f(V_\nu)\quad \Rightarrow \quad
s_\lambda s_\mu = c_{\lambda\mu}^{\;\;\nu}\,f(V_\nu)
$$
we can find the quantity $f(V_\lambda)$ for $\lambda$
corresponding to a given couple $(m,k)$ provided that all
$f(V_\mu)$ where $\mu$ corresponds to the couples $(l,r)$ such
that $l<m$ or $l=m,\,\,r<k$ are already known. \hfill
\rule{6.5pt}{6.5pt}

In what follows we will use this proposition for a proper
definition of the categorical trace. It will be defined as such a
morphism (\ref{tr}) that the corresponding dimension becomes an
\mbox{a-m} functional. By virtue of the above proposition we
should only find the quantities \r{val} for $\lambda=(1^k)$.

Describe now this construction in detail. Begin with the basic
object $V$. Since the space $V$ is finite dimensional, we can
identify  ${\rm End}(V)$ and $V\ot V^*$ in the usual way:
$$
\forall \, {\cal F}\in {\rm End}(V)\quad \longleftrightarrow\quad
e_i\ot F^i_{\,j}e^j\in V\ot V^*
$$
$F^i_{\,j}$ being a matrix of $\cal F$ in the basis $\{e_i\}$.
Define the mapping ${\rm tr}_V:\; {\rm End}(V) \rightarrow \Bbb K$
as the following composition of morphisms: \be {\rm tr}_V =
\alpha\,\ev{V}\circ \bar R_{V,V^*}.\label{def:tr-V} \ee

Note, that ${\rm tr}_V$ is defined by the above relation up to an
arbitrary nonzero factor $\alpha$ which will be specified later
from the requirement that the corresponding dimension would be an
a-m functional. In the fixed basis $\{e_i\}$ of $V$ one can write
this mapping in an explicit form. Taking $\bar R_{V,V^*}$ from
Proposition~\ref{star-braid} we find: \be \forall {\cal F}\in
{\rm End}(V):\quad e_i\rightarrow e_jF^j_{\,i}\quad
\Rightarrow\quad {\rm tr}_V({\cal F}) = \alpha\, {\rm Tr}(F\cdot
C) \label{expl-tr} \ee where the symbol $\rm Tr$ means the usual
matrix trace and the $n\times n$ matrix $C$ is defined to be
$$
C^i_{\,j}:= \sum_{a=1}^n Q^{ia}_{\;ja},\qquad Q\equiv (\pm
q)^{-\frac{1}{p}}\bar Q.
$$
Here the matrix $\bar Q$ is taken from (\ref{R-Q}) and the factor
$(\pm q)^\frac{1}{p}$ is included into $\alpha$.

The matrix $C$ has the following useful properties \cite{GPS}:
\begin{eqnarray}
&&R_{12}C_1C_2 = C_1C_2R_{12}\label{RCC}\\
&& C^{|1\rangle}_{\,\langle1'|} =
\frac{p_q}{q^p}\,u_{\langle1'23\dots p|}v^{|123\dots
p\rangle}\quad\Rightarrow\quad
{\rm Tr}\,C = \frac{p_q}{q^p}.\label{Tr-C}\\
&&{\rm Tr}_{(2)}R_{12}C_2 = {\rm id}_{(1)}\label{tr-2-R}.
\end{eqnarray}

The trace ${\rm tr}_{V^{\ot k}}$ in ${\rm End}(V^{\ot k})$ is
defined (also up to an arbitrary factor $\beta_{(k)}$) as the
mapping $V^{\ot k}\ot V^{*\ot k}\rightarrow \Bbb K$: \be {\rm
tr}_{V^{\ot k}} = \beta_{(k)}\,{\rm tr}_V^{(1)}\circ{\rm
tr}_V^{(2)}\dots\circ{\rm tr}_V^{(k)}\label{new-tr} \ee where
${\rm tr}_V^{(i)}:\,V^{\ot i}\ot V^{*\ot i}\rightarrow V^{\ot
(i-1)}\ot V^{*\ot (i-1)}$ reads
$$
{\rm tr}_V^{(i)} = {\rm id}^{\ot (i-1)}\ot \alpha^{-1}\,{\rm
tr}_V\ot {\rm id}^{\ot (i-1)}.
$$
For any basic object $V_\lambda$ $\lambda\vdash m$ one can find
the explicit form of the trace in ${\rm End}(V_\lambda)\cong {\rm
End}(V_\lambda(i))$ using definitions (\ref{space}),
(\ref{space*}) and property (\ref{RCC}): \be \forall\,{\cal
F}\in{\rm End}(V_\lambda(i)):\quad {\rm tr}_{V_\lambda} ({\cal
F}) = \beta_{(m)}\, {\rm Tr}_{(1\dots m)}(F\cdot C_\lambda)
\label{new-tra} \ee where ${\rm Tr}_{(1\dots m)}$ is the the
usual trace in the matrix space ${\rm Mat}_{n\times n}^{\ot
m}(\Bbb K)$ and the matrix $C_\lambda$ is of the form: \be
C_\lambda = {\cal Y}^\lambda_{ii}(R_1\dots R_{m-1})C_1C_2\dots
C_m. \label{C-lam} \ee {\remark{At the first sight the trace
${\rm tr}_{V_\lambda}$ depends on a concrete embedding $V_\lambda
\rightarrow V_\lambda(i)$. But in fact our definition is
invariant w.r.t. different embeddings of $V_\lambda$ into $V^{\ot
m}$. Indeed, as was already noticed in the previous section, any
primitive idempotents $Y^\lambda_{ii}$ of the Hecke algebra $H_m$
can be transformed into another idempotent $Y^\lambda_{jj}$ (with
the same $\lambda$) with the help of an invertible element
$X(i,j|\lambda)$ from $H_m$. The image ${\cal X}
=\rho_R(X(i,j|\lambda))$ of the element $X(i,j|\lambda)$ under
the representation $\rho_R$ (\ref{Rep}) represents a first kind
morphism of SW(V) and according to (\ref{RCC}) the string of
matrices $C_1\dots C_m$ in (\ref{C-lam}) commute with ${\cal X}$.
Therefore trace (\ref{new-tra}) is actually independent of the
index $i$ labelling the concrete embedding $V_\lambda\rightarrow
V_\lambda(i)$.}}

Find now the conditions on arbitrary factors $\alpha$ and
$\beta_{(k)}$ which would guarantee the dimensions defined via the
above trace to be an a-m functional.

Let us note that although the dimension defined by the usual trace
is an a-m functional such a trace is not a morphism of the
category if the braiding differs from the flip\footnote{ For
example, if our category SW(V) is supplied with the comodule
structure over the $RTT$ Hopf algebra (\ref{RTT}) then the usual
trace is not compatible with the comodule structure.}. This is
the reason why we should modify the usual trace by means of the
matrix $C$ (and all its extensions $C_\lambda$) multiplied by a
proper factor $\alpha$ (\ref{expl-tr}).

The fact that the matrix $C$ extends to ${\rm End}(V^{\ot m})$ in
a "group like" way
$$
C\;\rightarrow \;C^{\ot m}
$$
ensures the conservation of the additive and multiplicative
structure of the functional \r{new-tr}  restricted on the
identity operator  (what is just categorical dimension) iff we put
\be \beta_{(k)} =\alpha^k.\label{beta} \ee

So, if we want the trace ${\rm tr}_{V_\lambda}$ to be a morphism
and giving rise to an a-m functional we have the only free
parameter $\alpha$ at our disposal. This parameter (a
normalization of the trace) can be found from the following
condition $\dim_qV_{(1^p)} = 1$ (see footnote~\ref{p-dim}). The
following lemma plays the key role in finding such a
normalization (the relation presented below was found in another
but equivalent form in \cite{H}).

{\lemma{\label{Omega}For any Hecke symmetry of the rank $p$ the
following relation holds: \be {\rm Tr}_{(1\dots k)}({\cal
A}^{(k)}C_1\dots C_k)= q^{-pk} \frac{p_q!}{k_q!(p-k)_q!},\quad
1\le k\le p\label{om} \ee where the matrix ${\cal A}^{(k)}$ is
the image of $k$-th order antisymmetrizer $A^{(k)}$ under the
representation $\rho_R$ (\ref{Rep}). }}
\smallskip

{\bf Proof.}\ \ For an arbitrary $m\in \Bbb Z$ introduce the
auxiliary notation
$$
{\cal R}_i(m):= R_i -\frac{q^m}{m_q}
$$
where $m_q$ is defined in (\ref{q-num}). With this notation the
matrix ${\cal A}^{(m)}(R_1,\dots ,R_{m-1})$ can be presented in
the form:
$$
{\cal A}^{(m)} = {\frac{(-1)}{m_q}}^{m-1}\!\!{\cal A}^{(m-1)}{\cal
R}_{m-1}(m-1){\cal R}_{m-2}(m-2)\dots {\cal R}_{1}(1), \qquad
{\cal A}^{(1)}:={\rm id}.
$$
For the sake of shortness we introduce the notation:
$$
{\rm Tr}_{q\,(1\dots k)}\,X:= {\rm Tr}_{(1\dots k)}\,(XC_1\dots
C_k).
$$
Now we prove (\ref{om}) by induction. For $k=1$ we get with the
help of (\ref{Tr-C})
$$
{\rm Tr}\,C = q^{-p}p_q
$$
which gives the base of induction. Suppose that the assertion of
the lemma is valid up to some integer $k<p$. Then using the simple
relations
\begin{eqnarray*}
&&{\rm Tr}_{q\,(k+1)}{\cal R}_{k}(k) =
-q^{-p}\,\frac{(p-k)_q}{k_q}\,{\rm
id}_{(k)}\\
&&{\cal A}^{(k)}{\cal R}_{k-1}(k-1)\dots {\cal R}_{1}(1) =
(-1)^{k-1}k_q\,{\cal A}^{(k)}.
\end{eqnarray*}
one can complete the induction
$$
\begin{array}{rcl}
{\rm Tr}_{q\,(1\dots k+1)}\,{\cal A}^{(k+1)}
&\hspace*{-2.5pt}=&\displaystyle \hspace*{-2.5pt}
\frac{(-1)^k}{(k+1)_q}{\rm Tr}_{q\,(1\dots k+1)}\left({\cal
A}^{(k)}{\cal
R}_{k}(k)\dots {\cal R}_1(1)\right)\\
\rule{0pt}{8mm}&\hspace*{-2.5pt}=&\hspace*{-2.5pt} \displaystyle
\frac{(-1)^k}{(k+1)_q}{\rm Tr}_{q\,(1\dots k)}\left({\cal
A}^{(k)}\left[ {\rm Tr}_{q\,(k+1)}{\cal R}_{k}(k)\right]{\cal
R}_{k-1}(k-1)\dots {\cal
R}_{1}(1)\right)\\
\rule{0pt}{8mm}&\hspace*{-2.5pt}=&\hspace*{-2.5pt}\displaystyle
(-1)^{k+1}\,q^{-p}\,\frac{(p-k)_q}{k_q(k+1)_q}\,{\rm
Tr}_{q\,(1\dots k)}\left(
{\cal A}^{(k)}{\cal R}_{k-1}(k-1)\dots {\cal R}_1(1)\right)\\
\rule{0pt}{8mm}&\hspace*{-2.5pt}=&\hspace*{-2.5pt}q^{-p}\displaystyle
\frac{(p-k)_q}{(k+1)_q}\;{\rm Tr}_{q\,(1\dots k)}{\cal A}^{(k)}  =
\frac{q^{-p(k+1)} \,p_q!}{(k+1)_q!(p-k-1)_q!}.\hfill
\rule{6.5pt}{6.5pt}
\end{array}
$$

This lemma implies that if we take $\alpha = q^p$ (which is
equivalent to simple renormalization $C\rightarrow q^pC$ in all
formulas for the trace) then we get $\dim_q V_{(1^p)}=1$ and the
categorical dimension related to the trace becomes an a-m
functional well defined on objects of our category and hence on
the corresponding Grothendieck semiring.

Let us summarize the above consideration in the following
proposition.

{\proposition{ For any basic object $V_\lambda$, $\lambda\vdash
m$, of the category SW(V) define the trace ${\rm tr}_{V_\lambda}$
as a morphism ${\rm End}(V_\lambda) \rightarrow \Bbb K$ of the
form:
$$
\forall\,{\cal F}\in {\rm End}(V_\lambda):\quad {\rm
tr}_{V_\lambda}({\cal F}) = q^{pm}\, {\rm Tr}_{(1\dots m)}(F\cdot
C_\lambda).
$$
Then such a trace is a morphism in the category SW(V) and the
corresponding categorical dimension is an a-m functional. The
dimension of any basic object $V_\lambda$ is as follows: \be {\rm
dim}_qV_\lambda :={\rm tr}_{V_\lambda}({\rm id }) =
s_\lambda(q^{p-1}, q^{p-3},\dots ,q^{1-p})\label{dim-V} \ee
$s_\lambda$ being the Schur symmetric function (polynomial) in $p$
variables. }}
\smallskip

{\bf Proof.}\ \  To prove (\ref{dim-V}) we use
Proposition~\ref{AM-func} and Lemma~\ref{Omega}. By virtue of the
trace definition one gets from (\ref{om}):
$$
{\dim}_qV_{(1^k)} = q^{pk}\,{\rm Tr}_{q\,(1\dots k)}{\cal A}^{(k)}
= \frac{p_q!}{k_q!(p-k)_q!} \equiv f_k .
$$
Therefore due to Proposition~\ref{AM-func} we should find the
roots of the polynomial $\phi(t)$ (\ref{polynom}) with the above
$f_k$. But as is well known the generating function for such $f_k$
(which are $q$-binomial coefficients) is as follows
$$
E_q(t):=\prod_{k=0}^{p-1}(q^{2k+1-p}+t)\equiv \phi(t).
$$
Therefore the roots of the polynomial $\phi(t)$ with the
coefficients $f_k = \dim_qV_{(1^k)}$ are the numbers
$(-q^{2k+1-p})$, $k=0,1,\dots,p-1$ and result (\ref{dim-V})
follows now from Proposition~\ref{AM-func}. \hfill
\rule{6.5pt}{6.5pt}

{\remark{\label{concl} We want to complete the paper with the
following observation concerning the Koszul complexes considered
in \cite{G}. As was shown there, the Poincar\'e series
$P_{\pm}(t)$
$$
P_\pm(t):=\sum_l\dim \Lambda^l_\pm t^l
$$
of "symmetric" and "skewsymmetric" algebras $\Lambda_\pm(V)$ (see
(\ref{Sym-alg})) satisfy the relation \be P_+(t)\,P(-t)=1.
\label{clrel} \ee

We would like to point out that if we replace the usual
dimensions in formula \r{clrel} by categorical ones relation
\r{clrel} will be still valid. Moreover, the same is true if we
replace the dimensions by the values of any a-m functional. It is
not surprising since this fact reflects the well known relation
between elementary and complete symmetric functions. And any a-m
functional is nothing but a specialization of symmetric functions.

Thus, without constructing any deformed Koszul complex we can
obtain some numerical characteristics of quantum objects merely
replacing the usual dimensions by their $q$-analogs (if we
disregard the property of the quantum differential to be a
morphism of the category).}}


\begin{thebibliography}{DGMMM}

\bibitem[B]{B} Banica, {\em A reconstruction result for
R-matrix quantization of $SU(n)$}, QA/9806063.

\bibitem[BW]{BW} J.W.Barrett and B.W.Westbery, {\em
Spherical Categories}, Advances in Mathematics 143 (1999), pp
357-375.

\bibitem[CP]{CP} V.Chari and
A. Pressley, {\em A guide to Quantum Groups},  Cambrige University
Press, 1994.

\bibitem[DM]{DM} P.Deligne and J.S.Miln, {\em Tannakian Categories}
in: P.Deligne, J.S.Miln, A.Ogus, Shing Kuang-yen {\em Hodges
Cycles, Motives, and Shimura Varieties}, Lect. Notes in Math. 900,
Berlin, Heidelberg, NY: Springer-Verlag, 1982, pp. 101-228.

\bibitem[DJ]{DipJ} R. Dipper and G.James, {\it Representations of
Hecke algebras of General Linear Groups.}
Proc. London Math. Soc., 52(3) (1986), pp 20-52.\\
R. Dipper and G.James, {\it Block and Idempotents of Hecke
algebras of General Linear Groups}, Proc. London Math.
Soc., 54(3) (1987), pp 57-82.

\bibitem[DP]{DP} A.Dold and D.Puppe, {\em Duality, trace and transfer}
(Russian), Topology (Moscow 1979), Trudy Mat. Ins. Steklov 154
(154), pp 81-97.

\bibitem[DL]{DL} M.Dubois-Violette and G.Launer, {\em The
quantum group of a nondegenerated bilinear form},
Phys.Lett. B245 (1990), pp 175-177.

\bibitem[FRT]{FRT} L.D.Faddeev, N.Yu.Reshetikhin and L.A.Takhtajan,
{\it Algebra i Analiz} 1 (1989) 178 -- 206 (in Russian);
English translation in: {\it Leningrad Math. J.} 1 (1990),
pp 193--226.

\bibitem[Go]{Go} M.D.Gould, {\em Quantum Groups and Diagonalization
of the Braid Generator}, Lett. Math. Phys. 24 (1992) pp. 183 --
196.

\bibitem[G]{G} D.Gurevich, {\em Algebraic aspects of the quantum
Yang-Baxter equation}, Leningrad Math. J. 2 (1991), pp. 801--828.


\bibitem[GM]{GM} D.Gurevich and Z.Mriss, {\em Schur-Weyl Categories and
Non-quasiclassical Weyl Type Formula}, Hopf algebras and
quantum groups (Brussels, 1998) Lec. Notes. in Pure and
Appl. Math. 209, Dekker, NY 2000, pp 131--158.

\bibitem[GPS]{GPS} D.Gurevich, P.Pyatov and P.Saponov, {\em Hecke
Symmetries and Characteristic Relations on Reflection
Equation Algebra}, Lett. Mat. Phys. 41 (1997) pp. 255--264.

\bibitem[HIOPT]{5avt} L.K.Hadjiivanov, A.P.Isaev, O.V.Ogievetsky,
P.N.Pyatov and I.T.Todorov, {\em Hecke algebraic properties of
dynamical $R$-matrices. Application to related quantum matrix
algebras}, {\sf q-alg/9712026}

\bibitem[H]{H} Phung Ho Hai {\em On matrix quantum groups of type
$A_n$}, Int. J. Math. 11 (2000), pp. 1115--1146.

\bibitem[KW]{KW} D.Kazhdan, H.Wenzl, {Reconstructing Monoidal Categories},
Adv. in Soviet. Math. 16, part 2 (1993), pp. 111--136.


\bibitem[M]{M} I.G.Macdonald, {\em Symmetric functions and Hall
Polynomials}, Claredon Press, Oxford, 1979.

\bibitem[McL1]{ML} S.MacLane, {\em Natural associativity and
commutativity}, Rice Univ. Stud. 49 (1963), pp. 28--46.

\bibitem[McL2]{MacL} S.MacLane, {\em Categories for the Working
Mathematician}, Graduate Texts in Mathematics, Springer Verlag,
Berlin 1971.

\bibitem[OP]{PyOg} O.Ogievetsky and P.Pyatov, {\it Lecture on Hecke
algebra}, based on lectures at the International School
"Symmetries and Integrable systems", Dubna, 8-11 June, 1999. JINR
Publ. Dept., 2000.

\bibitem[T]{T} V.G.Turaev, {\em Quantum invariants of Knots and
3-Manifolds},  Walter de Gruyter, Berlin, NY 1994.

\end{thebibliography}
\end{document}